\documentclass[10pt,journal,compsoc]{IEEEtran}
\ifCLASSINFOpdf
\else
\fi
\usepackage{amsmath}
\usepackage{amssymb}
\usepackage{mathrsfs}
\usepackage{graphicx}
\usepackage{subfigure,bm}
\usepackage{amsmath}
\usepackage{amssymb}
\usepackage{algorithm}
\usepackage[switch]{lineno}
\usepackage{booktabs} 
\usepackage{algpseudocode}
\usepackage{color}

\newcommand{\EE}{\mathbb{E}}
\newcommand{\PP}{\textbf{Proj}}

 \def\biggiven{{\,\Big|\,}}

\newcommand{\Ac}{\mathcal{A}}
\newcommand{\Sc}{\mathcal{S}}
\newcommand{\Rc}{\mathcal{R}}
\newcommand{\Pc}{\mathcal{P}}
\newcommand{\Tc}{\mathcal{T}}
\newcommand{\RR}{\mathbb{R}}
\def \ttheta {{\boldsymbol \theta}}

\newtheorem{lemma}{Lemma}

\newtheorem{theorem}{Theorem}

\newtheorem{assumption}{Assumption}
\newtheorem{proposition}{Proposition}
\begin{document}

\title{  Adaptive Temporal Difference Learning with\\ Linear Function Approximation }
%
%
%

\author{
  Tao Sun$^*$\\
  College  of Computer\\
  National University of Defense Technology\\
  Changsha, Hunan 410073, China \\
  \texttt{nudtsuntao@163.com} \\
  \And Han Shen$^*$\\
   Department of ECSE\\
 Rensselaer Polytechnic Institute,  Troy \\
  Troy, NY, USA \\
  \texttt{shenh5@rpi.edu} \\
    \And Tianyi Chen\\
   Department of ECSE\\
 Rensselaer Polytechnic Institute,  Troy \\
  Troy, NY, USA \\
  \texttt{chent18@rpi.edu} \\
  \And
  Dongsheng Li\\
   College  of Computer\\
  University of California, Los Angeles \\
  Changsha, Hunan 410073, China\\
  \texttt{dsli@nudt.edu.cn} \\
}

\author{Tao Sun, Han Shen, Tianyi Chen, and
        Dongsheng Li
\thanks{The work of T. Sun and D. Li is sponsored in part by the National Key R\&D Program of China under Grant (2018YFB0204300)  and the National Natural Science Foundation of China under Grants (61932001 and 61906200).

Tao Sun and Dongsheng Li are with the
College of Computer, National University of Defense Technology,
Changsha, 410073, Hunan,  China (e-mails: \texttt{nudtsuntao@163.com},  \texttt{dsli@nudt.edu.cn}).

Han Shen and Tianyi Chen are with    Department of ECSE,
 Rensselaer Polytechnic Institute,
  Troy, NY, USA,
  (e-mails: \texttt{shenh5@rpi.edu,chent18@rpi.edu}).

 The first and second authors contributed equally.

Dongsheng Li is the corresponding author.}
}

%
%

\markboth{Journal of \LaTeX\ Class Files,~Vol.~14, No.~8, August~2015}%
{Shell \MakeLowercase{\textit{et al.}}: Bare Demo of IEEEtran.cls for IEEE Journals}
%



\maketitle

\begin{abstract}
This paper revisits the temporal difference (TD) learning algorithm for the policy evaluation tasks in reinforcement learning.
Typically, the performance of  TD(0) and TD($\lambda$) is very sensitive to the choice of stepsizes. Oftentimes, TD(0) suffers from slow convergence.
Motivated by the tight link between the TD(0) learning algorithm and the stochastic gradient methods, we develop a  provably convergent  adaptive  projected variant of the TD(0) learning algorithm with linear function approximation that we term AdaTD(0). In contrast to the TD(0), AdaTD(0) is robust or less sensitive to the choice of stepsizes. Analytically, we establish that to reach an $\epsilon$ accuracy, the number of iterations needed is $\tilde{O}(\epsilon^{-2}\ln^4\frac{1}{\epsilon}/\ln^4\frac{1}{\rho})$   in the general case, where $\rho$ represents the speed of the underlying Markov chain converges to the stationary distribution. This implies that the iteration complexity of AdaTD(0) is no worse than that of TD(0) in the worst case.   When the stochastic semi-gradients are sparse, we provide theoretical acceleration of AdaTD(0).
Going beyond TD(0), we develop an adaptive variant of TD($\lambda$), which is referred to as AdaTD($\lambda$).  Empirically, we evaluate the performance of AdaTD(0) and AdaTD($\lambda$) on several standard reinforcement learning tasks, which demonstrate the effectiveness of our new approaches.
\end{abstract}

\begin{IEEEkeywords}
Temporal Difference, Linear Function Approximation,  Adaptive Step Size,  MDP,  Finite-Time Convergence
\end{IEEEkeywords}

\section{Introduction}
Reinforcement learning (RL) involves a sequential decision-making procedure, where an agent takes (possibly randomized) actions in a stochastic environment over a sequence of time steps, and aims to maximize the long-term cumulative rewards received from the interacting environment. Owing to its generality, RL has been widely studied in many areas, such as control theory, game theory, operations research, multi-agent systems \cite{sutton1998introduction}.
Temporal Difference (TD) learning is one of the most commonly used algorithms for policy evaluation in RL \cite{sutton1988learning}.
TD learning provides an iterative procedure to estimate the value function with respect to a given policy based on samples from a Markov chain.
The classical TD algorithm adopts a tabular representation for the value function, which stores value estimates on a per state basis.
In large-scale settings, the tabular-based TD learning algorithm can
become intractable due to the increased number of states, and thus function approximation techniques are often combined with TD for better scalability and efficiency \cite{baird1995,tsitsiklis1996analysis}.

The idea of TD learning with function approximation is essentially to parameterize the
value function with a linear or nonlinear combination of fixed basis functions induced by the states that are termed feature vectors, and estimates the combination parameters in the same spirit of the tabular TD learning.
Similar to all other parametric stochastic optimization algorithms, however, the performance of the TD learning algorithm with function approximation is very sensitive to the choice of stepsizes. Oftentimes, it suffers from slow convergence \cite{even2003learning}. Ad-hoc adaptive modification of TD with function approximation has  often been used empirically, but their convergence behavior and  rate have not been fully understood.  When implementing TD-learning,  many practitioners use the  adaptive optimizer but without theoretical guarantees. This paper is devoted to the development of a provably convergent adaptive algorithm to accelerate the  TD(0) and TD($\lambda$) algorithms.
The key difficulty here is that the update used in the original TD does not follow the (stochastic) gradient direction of any objective function in an optimization problem, which prevents the use of the popular gradient-based optimization machinery. And the Markovian sampling protocol naturally involved in the TD update further complicates the analysis of adaptive and accelerated optimization algorithms.

\subsection{Related works}
We first briefly review related works in both the areas of TD learning and adaptive stochastic gradient.

\textbf{Temporal difference learning.}
The great empirical success of TD \cite{sutton1988learning} motivated active studies on the theoretical foundation of TD.
The first convergence analysis of TD was given by \cite{jaakkola1994convergence} using stochastic approximation techniques.
 In \cite{tsitsiklis1996analysis}, the characterization of limit points in TD  with linear function approximation has been studied, giving new intuition about the dynamics of TD learning.
  The ODE-based method (e.g., \cite{borkar2000ode}) has  dramatically inspired the subsequent development of research on asymptotic convergence of TD.
  Early convergence results of TD learning were mostly asymptotic, e.g., \cite{sutton2009nips}, because the TD update does not follow the (stochastic) gradient direction of any fixed objective function.
Non-asymptotic analysis for the gradient TD --- a variant of the original TD has been first studied by \cite{liu2015finite}, in which
the authors
reformulate  the original problem as new primal-dual saddle point optimization.
The finite-time analysis of TD with linear function approximation under i.i.d observation has been studied in \cite{dalal2018finite}; in particular, it is assumed
 that  observations in each iteration of TD   are independently drawn  from the steady-state distribution. In a concurrent line of research, TD has been considered in the view of the stochastic linear system, whose improved results are given by \cite{lakshminarayanan2018linear}.
Even without any fixed objective function to optimize, the proofs of \cite{dalal2018finite,lakshminarayanan2018linear} still follow a Lyapunov analysis like the SGD  due to the i.i.d sampling assumption and quadratic functions structure.
Nevertheless,   a more realistic assumption  for the data sampling in TD is the Markov  rather than the i.i.d process.
The finite-time convergence analysis under Markov sampling  is first presented in \cite{bhandari2018finite},
whose results are based on  the controls of
the gradient basis [Lemma 9, \cite{bhandari2018finite}] and a coupling [Lemma 11, \cite{bhandari2018finite}].
  The finite-time analysis for stochastic linear system under the Markov sampling is established by \cite{srikant2019finite,hu2019nips}, which is a general formulate of TD and applies potentially to other   problems.
 The finite-time analysis of multi-agent TD is proved by \cite{doan2019convergence}. However, all the aforementioned work leverages the original TD update.
In \cite{devraj2017zap}, the authors proposed   an improvement of Q-learning, which can  be used to TD, but with only asymptotic analysis being provided.
An adaptive   variant  of two time-scale stochastic approximation  was introduced by \cite{gupta2019finite}, that can also be applied to TD.
In \cite{vieillard2020momentum}, the authors introduce the momentum techniques for reinforcement learning and
 an extension on DQN.  In \cite{shani2020adaptive}, the authors
propose an adaptive scaling mechanism for TRPO and show that
it is the ``RL version" of traditional trust-region methods from convex analysis.


\textbf{Adaptive stochastic gradient descent.}
In machine learning areas different but related to RL, adaptive stochastic gradient descent methods have been actively studied. The first adaptive gradient (AdaGrad) is proposed by \cite{duchi2011adaptive,mcmahan2010adaptive}, and the algorithm demonstrated impressive numerical results when the gradients are sparse. While the original AdaGrad has a performance guarantee only in the convex case, the nonconvex AdaGrad has been studied by \cite{li2018convergence}. Besides the convex results, sharp analysis for nonconvex AdaGrad has also been  investigated in \cite{ward2018adagrad}.
Variants of AdaGrad have been developed in \cite{tieleman2012lecture,zeiler2012adadelta}, which
use alternative updating schemes (the exponential moving average
schemes) rather than the average of  the square of the past iterate.
The momentum technique applied to the adaptive stochastic algorithms gives birth to Adam and Nadam \cite{kingma2014adam,dozat2016incorporating}.
However, in \cite{reddi2019convergence}, the authors demonstrate that Adam may diverge under certain circumstances, and provide a new convergent Adam algorithm called AMSGrad.
Another method given by \cite{zou2019sufficient} is the use of decreasing factors for moving the average of the square of the past iterates.
 In \cite{chen2018convergence},  the convergence for generic Adam-type algorithms has been studied, which contains various adaptive methods.

\subsection{Comparison with existing analysis}
Our analysis considers the Markov sampling setting and thus differs \cite{dalal2018finite,lakshminarayanan2018linear}. It is worth mentioning that
the stepsizes chosen in this paper are different from that of \cite{gupta2019finite}: in \cite{gupta2019finite}, the ``adaptive" means that
the learning rate is reduced by multiplying a preset factor when
transient error is
dominated by the steady-state error; while in our paper, the ``adaptive" inherits the notion from Ada training method, i.e., using the learning rate associated  with the past gradients.
Our analysis cannot directly follow the techniques given by \cite{bhandari2018finite,srikant2019finite,hu2019nips,devraj2017zap} since  the learning rate in our algorithm is statistically
dependent on past information, which breaks previous Lyapunov analysis.

Thus,
this paper  uses
a delayed expectation technique  rather than bounding  a coupling [Lemma 11, \cite{bhandari2018finite}].
Furthermore, our analysis needs to deal with several terms related to adaptive style learning rates, which are quite complicated but absent in \cite{bhandari2018finite}.
Since the adaptive learning rate consists of past iterates  instead of being preset, these terms do not enjoy simple explicit bounds.
To this end, in this paper, we develop novel techniques (Lemmas \ref{lemmayu} and \ref{core1}).
On the other  hand, our analysis is also different from the adaptive SGD because TD update is not the stochastic gradient of any objective function; it uses biased samples generated from the Markov chain.

\subsection{Our contributions}
Complementary to existing theoretical RL efforts, we propose \emph{the first  provably convergent adaptive projected variant} of the TD learning algorithm with linear function approximation that has  finite-time convergence guarantees. For completeness of our analytical results, we investigate both the   TD(0) algorithm as well as the  TD($\lambda$) algorithm. In a nutshell, our contributions are summarized in threefold:

\textbf{c1)} We develop the adaptive variants of the TD(0) and TD($\lambda$) algorithms with linear function approximation. The new algorithms AdaTD(0) and AdaTD($\lambda$) are simple to use.

\textbf{c2)} We establish the finite-time convergence guarantees of AdaTD(0) and AdaTD($\lambda$), and they are not worse than those of TD and TD($\lambda$) algorithms in the worst case.

\textbf{c3)} We test our AdaTD(0) and AdaTD($\lambda$) on several standard RL benchmarks  and  show how these compare favorably to existing alternatives  like TD(0), TD($\lambda$), etc.

\section{Preliminaries}
This section introduces the notation, assumptions about the underlying MDP, and the setting of TD learning with linear function approximation.

\noindent \textbf{Notation:} The coordinate $j$ of a vector ${\bf x}$ is denoted by ${\bf x}_j$ and ${\bf x}^{\top}$ is transpose of ${\bf x}$. We use $\EE[\cdot]$ to denote the expectation with respect to the underlying probability space, and $\|\cdot\|$ for $\ell_2$ norm. Given a constant $R>0$ and ${\bf y}\in\RR^d$, $\PP_{R}({\bf y})$ denotes the projection of ${\bf y}$ to the ball $\{{\bf x}\in\RR^d\mid\|{\bf x}\|\leq R\}$.  For a matrix $A\in\RR^{S\times d}$, $\PP_{A}({\bf y})$ denotes the projection to space $\{A{\bf x}\mid {\bf x}\in\RR^d\}$.
 We denote the sub-algebra as $\sigma^k:=\sigma(\ttheta^0,\ttheta^1,\ldots,\ttheta^k)$, where $\ttheta^k$ is the $k$th iterate. We use $\tilde{O}(b)$ and $\tilde{\Theta}(b)$ to  hide the logarithmic factor of $b$.

\subsection{Markov Decision Process}
Consider a Markov decision process (MDP) described as a tuple  $(\Sc, \Ac, \Pc, \Rc, \gamma$), where $\Sc$ denotes the state space, $\Ac$ denotes the action space, $\Pc$ represents the transition matrix, $\Rc$ is the reward function, and $0<\gamma<1$ is the discount factor.  In this case, let $\mathcal{P}(s'|s)$ denote the transition probability from state $s$ to state $s'$. The corresponding transition reward is $\Rc(s,s')$.
We consider the finite-state case, i.e., $\Sc$ consists of  $|\mathcal{S}|$ elements, and a stochastic policy $\mu:{\cal S}\rightarrow {\cal A}$ that specifies an action given the current state $s$.
We use the following two assumptions on the stationary distribution and the reward.
\begin{assumption}\label{ass0}
  The transition  rewards are  uniformly bounded, i.e.,
$$
|\Rc(s,s')| \leq B,~\forall \,\, s,s' \in \mathcal{S}.
$$
\end{assumption}
\begin{assumption}\label{ass1}
  For any two states $s, s'\in \Sc$, it holds that
  $$\pi(s') = \lim_{t\to \infty} \Pc(s_t= s'|s_0=s)>0.$$
  There exist $\bar{\kappa}>0$ and $0\leq \rho<1$ such that
  $$\sup_{s\in\Sc}\{\sum_{s'\in\Sc}|\Pc(s_{t}=s'\mid s_0=s)-\pi(s')|\}\leq\bar{\kappa}\rho^t.$$
\end{assumption}
Assumptions \ref{ass0} and \ref{ass1} are standard in MDP. For irreducible and aperiodic Markov chains, Assumption \ref{ass1} can always hold \cite{levin2017markov}.  The constant $\rho$  represents the Markov chain's speed accessing the stationary distribution $\pi$.
When the number of states is finite, the Markovian transition kernel is a matrix $\mathcal{P}$, and $\rho$ is identical to the second largest eigenvalue of $\mathcal{P}$.
An important notion in the Markov chain is the mixing time, which measures the time that a Markov chain needs for its current state distribution roughly matches the stationary one $\pi$. Given an $\epsilon>0$, the mixing time is defined as  $$\tau(\epsilon):=\min_{t\in \mathbb{Z}^+}\{\bar{\kappa}\rho^t\leq \epsilon\}.$$ With Assumption \ref{ass1}, we can see $\tau(\epsilon)=O({\ln\frac{1}{\epsilon}}/{\ln\frac{1}{\rho}})$. That means if $\rho$ is small, the mixing time is short.

This paper considers the \emph{on policy} setting, where both the target  and  behavior policies are $\mu$.
For a given policy $\mu$, since the actions or the distribution of actions will be uniquely determined, we thus eliminate the dependence on the action in the rest of the paper.
We denote the expected reward at a given state $s$ by $$\Rc(s):= \sum_{s'\in \mathcal{S}}\Pc(s'|s) \Rc(s,s').$$
The value function $V_{\mu}:\Sc \to \mathbb{R}$ associated with a policy $\mu$ is the expected cumulative discounted reward from a given state $s$, that is
$$
V_\mu(s) = \mathbb{E}\left[ \sum_{t=0}^\infty \gamma^t \Rc(s_t) \biggiven s_0=s \right],
$$
where the expectation is taken over the trajectory of states generated under $\Pc$ and $\mu$.
The restriction on discount $0<\gamma<1$ can guarantee   the boundedness of  $V_\mu(s)$.
The Markovian property of MDP yields the well-known Bellman equation
\begin{equation}\label{eq.bell}
	\Tc_\mu V_\mu = V_\mu,
\end{equation}
where the  operator $\Tc_\mu$ on $V$ is defined as
$$
(\Tc_\mu V)(s):= \Rc(s)+ \gamma \sum_{s' \in \mathcal{S}} \Pc(s'|s) V(s'), \,\, s\in \Sc.	
$$

Solving the (linear) Bellman equation allows us to find the value function  $V_{\mu}$ induced by a given policy $\mu$. However, in practice, $|\mathcal{S}|$ is usually very large and
 computationally intractable.
 An alternative method is to leverage the linear \cite{sutton1998introduction} or non-linear  function approximations (e.g., kernels and neural networks \cite{mnih2015human}). We focus on the linear case here, that is
\begin{equation}
	V_{\mu}(s) \approx V_{\ttheta}(s):= \phi(s)^\top \ttheta,
\end{equation}
where $\phi(s) \in \mathbb{R}^d$ is the feature vector for state $s$, and $\ttheta \in \mathbb{R}^d$ is a parameter vector.  To reduce difficulty caused by the dimension, $d$ is set smaller than $|\mathcal{S}|$.
With the linear   function  approximator, the vector $V_{\ttheta}\in \mathbb{R}^{|\mathcal{S}|}$ becomes 
\begin{align*}
V_{\ttheta} = \Phi \ttheta,
\end{align*}
where the feature matrix is defined as
$$\Phi:= [\phi(s_1), \phi(s_2),\ldots, \phi(s_{|\mathcal{S}|})]^{\top}\in \mathbb{R}^{{|\mathcal{S}|} \times d} $$
with $s_n$ being the $n$th state.

\begin{assumption}\label{ass2}
  For any state $s\in\Sc$, we assume the feature vector is uniformly bounded such that
  $\|\phi(s)\|\leq 1$, and the feature matrix $\Phi $ is full column-rank.
\end{assumption}

It is not hard to guarantee Assumption \ref{ass2} since the feature map $\phi$ is chosen by the users and $|\mathcal{S}|>d$. With Assumptions \ref{ass1} and \ref{ass2}, we can see that the matrix $\Phi^{\top}\textrm{Diag}(\pi)\Phi$ is positive define, and we denote its minimal eigenvalue as follows
$$
	\omega:=\lambda_{\min}\left(\Phi^{\top}\textrm{Diag}(\pi)\Phi\right)>0.
$$
With the linear approximation of value function, the task then is tantamount to finding $\ttheta\in\RR^d$ that obeys the Bellman equation given by
$$
	\Phi\ttheta=\Tc_\mu \Phi\ttheta.
$$
However, $\ttheta$ that satisfies such an equation may not exist if $ V_{\mu}\notin \{\Phi\ttheta\mid\ttheta\in\RR^d\}$.
Instead, there always exists a unique solution $\ttheta^*$ for the projected Bellman equation \cite{tsitsiklis1996analysis}, given by
\begin{equation}
		\Phi\ttheta^*=\PP_{\Phi}(\Tc_\mu \Phi\ttheta^*),
\end{equation}
where $\PP_{\Phi}$ is the projection  onto the span of $\Phi$'s columns.

\section{Adaptive Temporal Difference Learning}

\subsection{TD with linear function approximation}
TD(0) algorithm starts with an initial parameter $\theta^0$. At iteration $k$, after sampling states $s_k$, $s_{k+1}$, and reward $r(s_{k+1},s_k)$ from a Markov chain, we can compute the TD (temporal difference) error which is also called the Bellman error:
\begin{equation}
    d_k := r(s_{k+1},s_k)+\gamma V_{\ttheta}(s_{k+1})-V_{\ttheta}(s_k).
\end{equation}
The TD error is subsequently used to compute the stochastic semi-gradient:
\begin{equation}\label{sampling}
    \overline{{\bf g}}(\ttheta^k;s_{k},s_{k+1}):=d_k\nabla V_{\ttheta^k}(s_k)=d_k\phi(s_{k}).
\end{equation}
The traditional TD(0) with linear function approximation performs SGD-like update as
\begin{equation}\label{eq.linear-TD}
		\ttheta^{k+1}=\ttheta^k+\eta \overline{{\bf g}}(\ttheta^k;s_{k},s_{k+1}).
\end{equation}
\begin{table}
      \begin{algorithm}[H]
    \caption{ Projected TD(0) Learning}\label{alg0}
	\begin{algorithmic}[1]
\State{\textbf{Parameters}: learning rate $\eta$.}
\State{\textbf{Initialization}: ${\bf g}^0=\textbf{0}$,  ${\bf m}^0=\textbf{0}$, $v^0=0$}
\For{$k=1,2,\ldots$}
\State{sample a state transition $s_{k}\rightarrow s_{k+1}$ from $\mu$}
\State{calculate $\overline{{\bf g}}(\ttheta^k;s_{k},s_{k+1})$ in \eqref{sampling}}
\State{update the parameter $\ttheta^k$ as \begin{equation}\label{eq.linear-projTD}
		\ttheta^{k+1}=\PP_{R}[\ttheta^k+\eta \overline{{\bf g}}(\ttheta^k;s_{k},s_{k+1})],
\end{equation}}
\EndFor
	\end{algorithmic}
  \end{algorithm}
\end{table}

The update  TD(0)  makes sense  because the direction $\overline{{\bf g}}(\ttheta;s_{k},s_{k+1})$ is a good one since it is asymptotically close to the direction whose limit point is $\ttheta^*$. Specifically, it has been established that \cite{tsitsiklis1996analysis}
$$
		\lim_{k\rightarrow\infty}\EE [\overline{{\bf g}}(\ttheta;s_{k},s_{k+1})]={\bf g}(\ttheta),
$$
%
where ${\bf g}(\ttheta)$ is defined as
$$
{\bf g}(\ttheta):=\Phi^{\top}\textrm{Diag}(\pi)(\Tc_{\mu}\Phi\ttheta-\Phi\ttheta).
$$
We  term ${\bf g}(\ttheta)$ as the limiting update direction,   ensuring that ${\bf g}(\ttheta^*)=\textbf{0}$. Note that while $\overline{{\bf g}}(\ttheta;s_{k},s_{k+1})$ is an unbiased estimate under the stationary $\pi$, it is not for a finite $k$ due to the Markovian property of $s_k$.
Therefore, the TD(0) update \eqref{eq.linear-TD} is   asymptotically akin to the  stochastic approximation
$$\ttheta^{k+1}=\ttheta^k+\eta \overline{{\bf g}}(\ttheta^k;s,s'),$$
where $s,s'$ are independently drawn from the stationary distribution $\pi$.

Nevertheless, an important property of the limiting direction ${\bf g}(\ttheta)$, found by \cite{tsitsiklis1996analysis}, is that: for any $\ttheta\in\RR^d$, we have
\begin{equation}\label{eq.TD-asymp}
	  \langle\ttheta^*-\ttheta,{\bf g}(\ttheta)\rangle\geq (1-\gamma)\omega\|\ttheta^*-\ttheta\|^2.
\end{equation}
 An important observation follows from this inequality readily: only one $\ttheta^*$ satisfies ${\bf g}(\ttheta^*)=\textbf{0}$. We can show this by contradiction.
If there exists another $\bar{\ttheta}$ such that ${\bf g}(\bar{\ttheta})=\textbf{0}$, we have
$0=\langle\ttheta^*-\bar{\ttheta},{\bf g}(\bar{\ttheta})\rangle\geq (1-\gamma)\omega\|\ttheta^*-\bar{\ttheta}\|^2$, which again means $\ttheta^*=\bar{\ttheta}$.

To ensure the boundedness of $\ttheta^k$ and simplify the convergence analysis, projection is used in \eqref{eq.linear-projTD} (see e.g., \cite{bhandari2018finite}).
In \cite{bhandari2018finite}, it has been shown that if $R\geq   2B/\sqrt{\omega}(1-\gamma)^{\frac{3}{2}}$ ($B$ has appeared  in Assumption 1), the projected TD(0) does not exclude all the limit points of the TD(0) (such a fact still holds for our proposed algorithms).
The finite-time convergence  of projected  TD(0) is analyzed by \cite{bhandari2018finite}. The projection step is removed in \cite{srikant2019finite}, with almost the same results being proved. But in \cite{srikant2019finite}, the authors just studied the constant stepsize, while \cite{bhandari2018finite} shows more cases, including diminishing stepsize cases. In this paper, we consider a more complicated scheme that cannot be analyzed using  techniques in \cite{srikant2019finite}. Thus, the projection is still needed.

\subsection{Adaptive TD development}
Motivated by the recent success of adaptive SGD methods such as \cite{duchi2011adaptive,tieleman2012lecture,reddi2019convergence}, this paper aims to develop an adaptive version of TD with linear function approximation that we term AdaTD(0) (adaptive TD).
Unlike TD(0) method, in which step size is often a constant, AdaTD(0) seeks to scale the step size depending on the norm of stochastic gradient. Additionally, instead of using a one-step gradient, AdaTD(0) utilizes momentum, which gives the algorithm a memory of history information.

As presented in the last section, $\overline{{\bf g}}(\ttheta;s_{k},s_{k+1})$ is a stochastic estimate of ${\bf g}(\ttheta)$. Based on this observation, we develop the adaptive scheme for TD(0).
Different from projected TD(0), we use the update direction ${\bf m}^k$, which is the exponentially weighted average of stochastic gradients. It is further scaled by $v^k$,  the moving average of the squared norm of stochastic semi-gradients. Intuitively, when the gradient is large, the algorithm will take smaller steps. To prevent a too large scaling, we use a positive hyper-parameter $\delta$ for numerical stability. The AdaTD(0) update is given by \eqref{adatd-full}.
The key difference between AdaTD(0) and the  TD(0) method is that AdaTD(0) utilizes the history information in the update of both first moment estimate ${\bf m}^k$ and second moment estimate $v^k$.
Unlike projected  TD(0) whose
   asymptotically expected TD update, ${\bf g}(\ttheta^k)$, is a good direction as is evident from \eqref{eq.TD-asymp},
the update direction ${\bf m}^k$ in AdaTD(0) is an exponentially weighted version of $\overline{{\bf g}}(\ttheta^k;s_{k},s_{k+1})$, which makes our analysis more challenging.

Although the variance term in AdaTD(0) is used as a sum form, it can be   rewritten  as an  exponentially weighted moving average. If we denote $\hat{ v}^k:=\frac{\sum_{i=1}^k \|{\bf g}^i\|^2}{k}=\frac{v^k}{k}$, the last two steps of Ada-TD(0) can be reformulated as
  	\begin{align*}
 \hat{ v}^k&=(1-\frac{1}{k})\hat{ v}^{k-1}+\frac{1}{k}\hat{ v}^{k},\\
  \ttheta^{k+1}& = \PP_{R}(\ttheta^k +\frac{\eta}{\sqrt{k}} {\bf m}^k/(\hat{v}^k+\frac{\delta}{k})^{1/2}),
\end{align*}
Note that in this form,  weights are   $\{1/k\}_{k\geq 1}$ and stepsizes are $\{\frac{\eta}{\sqrt{k}}\}_{k\geq 1}$, which obey the sufficient conditions to guarantee the convergence of  Adam-type algorithms with exponentially weighted average variance term \cite{zou2019sufficient,chen2018convergence}.

\begin{table*}
\vspace{-0.6cm}
    \begin{tabular}{c c}
    \hspace{-0.1cm}
\begin{minipage}[t]{8.5cm}
  \vspace{0pt}
  \begin{algorithm}[H]
    \caption{ Projected Adaptive TD(0) Learning}\label{alg1}
	\begin{algorithmic}[1]
\State{\textbf{Parameters}: $\eta, \delta$,  $\beta$, $\gamma$, $R$.}
\State{\textbf{Initialization}: ${\bf g}^0=\textbf{0}$,  ${\bf m}^0=\textbf{0}$, $v^0=0$}
\For{$k=1,2,\ldots$}
\State{sample a state transition $s_{k}\rightarrow s_{k+1}$ from $\mu$}
\State{calculate $\overline{{\bf g}}(\ttheta^k;s_{k},s_{k+1})$ in \eqref{sampling}}
\State{update the parameter $\ttheta^k$ as \begin{subequations}\label{adatd-full}
	\begin{align}
    {\bf m}^k &=\beta {\bf m}^{k-1} + (1-\beta) \overline{{\bf g}}(\ttheta^k;s_{k},s_{k+1}),\\
  v^k &=   v^{k-1} +  \|\overline{{\bf g}}(\ttheta^k;s_{k},s_{k+1})\|^2,\\
  \ttheta^{k+1}& = \PP_{R}(\ttheta^k +\eta {\bf m}^k/\sqrt{v^k+\delta}).\label{adatd-scheme}
\end{align}
\end{subequations}}
\EndFor
	\end{algorithmic}
  \end{algorithm}
\end{minipage}
&
\begin{minipage}[t]{8.5cm}
  \vspace{0pt}
  \begin{algorithm}[H]
    \caption{ Projected Adaptive TD($\lambda$) Learning}\label{alg2}
	\begin{algorithmic}[1]
\State{\textbf{parameters}: $\eta$,  $\beta$, $\gamma$, $\delta$, $\lambda$, $\hat{R}$.}
\State{\textbf{initialization}: ${\bf g}^0=\textbf{0}$, ${\bf z}^0=\textbf{0}$, ${\bf m}^0=\textbf{0}$, $v^0=0$}
\For{$k=1,2,\ldots$}
\State{sample a state transition $s_{k}\rightarrow s_{k+1}$ from $\mu$}
\State{update ${\bf z}^k=(\gamma\lambda){\bf z}^{k-1}+\phi(s_k)$}
\State{update \begin{align}\label{sampling2}
        &\overline{{\bf g}}^\lambda(\ttheta;s_{k},s_{k+1},{\bf z}^k):=r(s_{k+1},s_k){\bf z}^k\nonumber\\
        &\qquad+\gamma\phi(s_{k+1})^{\top}\ttheta {\bf z}^k-\phi(s_{k})^{\top}\ttheta{\bf z}^k\in \RR^d,
\end{align}}
\State{update the parameter $\ttheta^k$ as \eqref{adatd-full} with $R\leftarrow\hat{R}$}
\EndFor
	\end{algorithmic}
  \end{algorithm}
\end{minipage}
   \end{tabular}

\end{table*}

\subsection{Finite-time analysis of projected AdaTD(0)}\label{section:analysis_proj_atd}
Because  the main results depend on constants related to the bounds, we present them
in Lemma \ref{lebound}.
\begin{lemma}\label{lebound}
The following bounds hold for $(\ttheta^k)_{k\geq 0}$ generated by AdaTD(0)
\begin{align}\label{bound1}
\|\ttheta^{k}-\ttheta^{*}\|\leq 2R,~~~\|{\bf g}^k\|  \leq G,~~~\|{\bf m}^k\|  \leq G
\end{align}
where we define ${\bf g}^k:=\overline{{\bf g}}(\ttheta^k;s_{k},s_{k+1})$, and $G:=2R+B$.
\end{lemma}
Lemma \ref{lebound} follows readily. The bounds presented in Lemma \ref{lebound} are critical for the subsequent analysis.

\medskip

The convergence analysis of AdaTD(0) is more challenging than that of both TD and adaptive SGD.
Compared with the analysis of adaptive SGD methods in e.g., \cite{duchi2011adaptive,tieleman2012lecture,reddi2019convergence}, even under the i.i.d. sampling, the stochastic direction $\overline{{\bf g}}(\ttheta^k;s_{k},s_{k+1})$ used in \eqref{adatd-full} fails to be the stochastic gradient of any objective function,
aside from the fact that samples are drawn from a Markov chain. Compared with TD(0), the actual update of $\ttheta^k$ in AdaTD(0) involves the history information of both the first and the second moments of $\overline{{\bf g}}(\ttheta^k;s_{k},s_{k+1})$, which makes the analysis of TD in e.g., \cite{bhandari2018finite,srikant2019finite} intractable.
\begin{theorem}\label{th1}
Suppose $(\ttheta^k)_{k\geq 0}$ are generated by AdaTD(0) with
$R\geq2B/\sqrt{\omega}(1-\gamma)^{\frac{3}{2}}$
under the Markovian observation.
Given   $\eta>0, \delta>0, 0\leq \beta<1$,   we have
\begin{align*}
&\min_{1\leq k\leq K}\EE(\|\ttheta^*-\ttheta^k\|^2)\\
&\qquad\leq \Big[C_1\ln\big(\frac{\delta+KG^2}{\delta}\big)\Big]/\sqrt{K}+C_2/\sqrt{K},
\end{align*}
where $C_1$ and  $C_2$  are given as
\begin{align*}
 C_1&:=\frac{16(\ln K/\ln\frac{1}{\rho})^2G}{\sqrt{\delta}(1-\gamma)^2\omega^2}+\frac{2\eta\beta G}{(1-\beta)(1-\gamma)\omega}\\
 &\quad+\frac{\eta G}{(1-\gamma)\omega}+\frac{4RG^2}{(1-\gamma)\omega\delta}=O(\frac{(\ln K/\ln\frac{1}{\rho})^2}{\sqrt{\delta}}),\\
   C_2&:=\frac{4R^2G}{(1-\gamma)\omega\eta}+\frac{4R\eta G^2}{\delta^{\frac{1}{2}}(1-\beta)(1-\gamma)\omega}\\
  &\quad+\frac{4R \bar{\kappa}(B+R\gamma+R)G }{ \sqrt{\delta}(1-\gamma)\omega}=O(\frac{1}{\sqrt{\delta}}).
\end{align*}
\end{theorem}
With Theorem 1,  to achieve $\epsilon$-accuracy for $\min_{1\leq k\leq K}\{\EE(\|\ttheta^*-\ttheta^k\|^2)\} $,  we need
\begin{align}\label{rela}
C_2/\sqrt{K}\leq \epsilon/2,~~~~~~
       \frac{C_1\ln\big(\frac{\delta+KG^2}{\delta}\big)}{\sqrt{K}}\leq  \epsilon/2.
\end{align}
Since $C_2= O(1)$, with the first inequality of \eqref{rela}, $K=\tilde{O}(\frac{1}{\epsilon^2})$.
Further, due to $C_1= O([\ln K/\ln\frac{1}{\rho}]^2)=O([\ln \frac{1}{\epsilon}/\ln\frac{1}{\rho}]^2)$, with the second inequality of \eqref{rela}, we have
$C_1\ln K/\sqrt{K}=\tilde{O}(C_1/\sqrt{K})=O(\epsilon)$.
Therefore, we obtain a solution whose square distance to $\ttheta^*$ is $\epsilon$, the iteration needed is
\begin{align}\label{complexity-ori}
		\tilde{O}\left(C_1^2/\epsilon^2\right)=\tilde{O}\left(\ln^4\frac{1}{\epsilon}\Big/(\epsilon^2\ln^4\frac{1}{\rho})\right).
\end{align}
When $\rho$ is much smaller than $1$, the term $\ln\frac{1}{\epsilon}/\ln\frac{1}{\rho}$ keeps at a relatively small  level. Recall that the state-of-the-art convergence result of TD(0) given in \cite{bhandari2018finite} is $\tilde{O}( \ln^2\frac{1}{\epsilon}\big/(\epsilon^2\ln^2\frac{1}{\rho}))$, the  rate of AdaTD(0) is close to TD(0) in the general case.
We do not present a faster speed for technical reasons. In fact, this is also the case for the adaptive SGD \cite{duchi2011adaptive,tieleman2012lecture,reddi2019convergence}.
Specifically, although numerical results demonstrate the advantage of adaptive methods, the worst-case convergence rate of adaptive methods is still similar to that of the stochastic gradient descent method.
To reach a desired error as $\min_{1\leq k\leq K}\{\EE(\|\nabla f({\bf x}^k)\|^2)\}\leq \epsilon$ in adaptive SGD, where ${\bf x}^k$ is the $k$th iterative and $\epsilon>0$, the iteration number $K$ needs to be set as $O(\frac{1}{\epsilon^2})$, which is identical to the SGD.

\textbf{Sketch of the proofs:}
Now, we present the sketch of the proofs for the main result. Because AdaTD(0) does not have any objective function to optimize, we consider sequence the $(\|\ttheta^{k}-\ttheta^*\|^2)_{k\geq 0}$.

Using the update \eqref{adatd-full}, we have
\begin{align}\label{division-pre}
    	&\|\ttheta^*-\ttheta^{k+1}\|^2-\|\ttheta^{*}-\ttheta^{k}\|^2\nonumber\\
    	&\leq  \frac{2\eta \langle{\bf m}^k,\ttheta^{k}-\ttheta^{*}\rangle}{(v^k+\delta)^{\frac{1}{2}}}
    	+\frac{\eta^2\|{\bf m}^k\|^2}{( v^k+\delta)}\nonumber\\	
    	&\leq   \underbrace{2\eta \langle \ttheta^{k}-\ttheta^{*},  {\bf m}^k\rangle/( v^{k-1}+\delta)^{\frac{1}{2}}}_{I_1^k} +\underbrace{\eta^2\|{\bf m}^k\|^2/( v^k+\delta)}_{I_2^k}\nonumber\\
& \quad+\underbrace{2\eta \langle \ttheta^{k}-\ttheta^{*},  {\bf m}^k\rangle  (1/(v^{k}+\delta)^{\frac{1}{2}}-1/( v^{k-1}+\delta)^{\frac{1}{2}})}_{I_3^k}.
\end{align}
Here, we split the term $\langle{\bf m}^k,\ttheta^{k}-\ttheta^{*}\rangle/(  v^k+\delta)^{\frac{1}{2}}$ as $I_1^k+I_3^k$ because both ${\bf m}^k$ and $v^k$ statistically depends on $\overline{{\bf g}}(\ttheta^k;s_{k},s_{k+1})$, which makes it hard to calculate the expectation.
 Subtracting $-I_1^k$ to both sides of   \eqref{division-pre}, we get
\begin{align}\label{division}
    	-I_1^k\leq  \|\ttheta^{*}-\ttheta^{k}\|^2-\|\ttheta^*-\ttheta^{k+1}\|^2+ I_2^k+I_3^k.
\end{align}
With  division \eqref{division}, the analysis can be divided into three parts: {\bf A)} determine the  lower bound of $-\EE [I_1^k]$ (i.e, the upper bound of $\EE [I_1^k]$);  {\bf B)}  prove the   summability of $(\EE [I_2^k])_{k\geq 1}$, i.e., $\sum_{k=1}^{+\infty}I_2^k<+\infty$; {\bf C)}  prove the  summability of $(\EE[I_3^k])_{k\geq 1}$.

{\bf A)} Calculating the  upper bound of $\EE {[I_1^k]}$ is the most difficult part, which consists of two substeps.
Since ${\bf m}^k$ is a convex combination of ${\bf m}^{k-1}$ and $\overline{{\bf g}}(\ttheta^k;s_{k},s_{k+1})$.
We recursively bound $\EE [I_1^k]$ by analyzing $\EE[\langle \ttheta^{k}-\ttheta^{*},  \overline{{\bf g}}(\ttheta^k;s_{k},s_{k+1})\rangle/( v^{k-1}+\delta)^{\frac{1}{2}}]$.
Assume that the Markov chain has reached its stationary distribution, in which, $\EE[\overline{{\bf g}}(\ttheta^k;s_{k},s_{k+1})]={\bf g}$.
 But the stationary distribution is not assumed in our setting.
Therefore, we need to
bound the difference between $\EE\overline{{\bf g}}(\ttheta^k;s_{k},s_{k+1})$ and ${\bf g}$. Unlike i.i.d setting, directly bounding them is difficult in the Markovian setting since it will lead to some uncontrollable error in the final convergence result. To solve this technical issue, we consider $\EE[\overline{{\bf g}}(\ttheta^{k-T};s_{k},s_{k+1})\mid\sigma^k]$ and ${\bf g}(\ttheta^{k-T})$. This is because although $s_k$ is biased, $\Pc(s_k|s_{k-T}=s)$ is very close to $\pi(s)$ when $T$ is large.
Using this technique, we prove the  following Lemma.
\begin{lemma}\label{legeo}
Assume $(\ttheta^k)_{k\geq 0}$ are generated by AdaTD(0). Given an integer $K_0\in \mathbb{Z}^+$, we have
\begin{align}
\left\|\EE\left[ \overline{{\bf g}}(\ttheta^{k-K_0};s_{k},s_{k+1})\right]-{\bf g}(\ttheta^{k-K_0})\right\|\leq  \kappa\rho^{K_0},
\end{align}
where the constant is defined as $\kappa:=\bar{\kappa}(B+R\gamma+R)$.
\end{lemma}
 In the second substep, because we have got the biased bound caused by the Markovian stochastic process, it is possible to bound the difference $\EE[\langle  \ttheta^{k}-\ttheta^{*},  \overline{{\bf g}}(\ttheta^k;s_{k},s_{k+1})\rangle/( v^{k-1}+\delta)^{\frac{1}{2}}]$ between $\EE[\langle  \ttheta^{k}-\ttheta^{*},  {\bf g}(\ttheta^k)\rangle/( v^{k-1}+\delta)^{\frac{1}{2}}]$ with the  following lemma.
\begin{lemma}\label{lemmayu}
  Given $K_0\in \mathbb{Z}^+$, we have
\begin{align}
&\EE\Big[\langle  \ttheta^{k}-\ttheta^{*}, \overline{{\bf g}}(\ttheta^k;s_{k},s_{k+1})\rangle/( v^{k-1}+\delta)^{\frac{1}{2}}\Big]\nonumber\\
&\leq \frac{1}{2}\EE\Big[\langle\ttheta^{k}-\ttheta^{*}, {\bf g}(\ttheta^{k})\rangle/( v^{k-1}+\delta)^{\frac{1}{2}}\Big]+\frac{2R\kappa \rho^{K_0}}{\sqrt{\delta}}\nonumber\\
& +\frac{8K_0}{\delta^{1/2}(1-\gamma)\omega}\sum_{h=K_0}^1\EE \Big[\|{\bf m}^{k-h}\|^2/(v^{k-h}+\delta)^{\frac{1}{2}}\Big].
\end{align}
\end{lemma}

{\bf B)} The second part is to bound $\sum_{k}\EE[\|{\bf m}^k\|^2/( v^k+\delta)]$. Because $v^k$ depends on $(\|{\bf g}^k\|^2)_{k\geq 0}$ (${\bf g}^k:=\overline{{\bf g}}(\ttheta^k;s_{k},s_{k+1})$), we expand $\|{\bf m}^k\|^2$ by $(\|{\bf g}^j\|^2)_{j\geq 0}$,
 i.e., \eqref{core0-t1} in Lemma \ref{core0} in the Appendix.
We then apply a provable result (Lemma  \ref{le1} in the Appendix) to the right side of \eqref{core0-t1} and get \eqref{core0-t2}.

\begin{lemma}\label{core0}
Let $(\hat{\bf m}_k)_{k\geq 0}$ be defined in \eqref{notation1}, we have
\begin{align}\label{core0-t1}
\sum_{k=1}^K\hat{\bf m}_k\leq   \sum_{j=1}^{K-1}\EE\Big[\|{\bf g}^j\|^2/( v^{j}+\delta)\Big].
\end{align}
Further, with the boundedness of $(\|{\bf g}^k\|)_{k\geq 0}$, we then get
\begin{align}\label{core0-t2}
\sum_{k=1}^K\hat{\bf m}_k\leq \ln(\frac{\delta+(K-1)G^2}{\delta}).
\end{align}
\end{lemma}

{\bf C)} While the third part is the easiest and obvious. The boundedness of the points indicates the uniform bound of  $(|\langle \ttheta^{k}-\ttheta^{*},  {\bf m}^k\rangle|)_{k\geq 0}$. The monotonicity of $(v^k)_{k}$ then yields the summable bound.

 Summing \eqref{division} from $k=1$ to $K$,
with the proved bounds in these three parts, we then can bound $\sum_{k=1}^K  \langle\ttheta^*-\ttheta^k,{\bf g}(\ttheta^k)\rangle/( v^{k}+\delta)^{\frac{1}{2}}$. Once with the descent property $\langle\ttheta^*-\ttheta,{\bf g}(\ttheta)\rangle\geq (1-\gamma)\omega\|\ttheta^*-\ttheta\|^2$, we can derive the main convergence result.

\medskip

 The acceleration of adaptive SGD is proved under extra assumptions like the sparsity of the stochastic gradients. In the following, we present the acceleration result of AdaTD(0) also with an extra assumption.

\begin{theorem}\label{th1a}
Suppose the conditions of Theorem \ref{th1} hold and
\begin{align}\label{fastass}
v^{k}\leq c k^{\nu},
\end{align}
where $c>0$ is a universal constant and $0<\nu\leq 1$.
 To achieve $\epsilon$-accuracy for $\min_{1\leq k\leq K}\{\EE(\|\ttheta^*-\ttheta^k\|^2)\} $,  the needed iteration is $\tilde{O}\left((\ln\frac{1}{\epsilon})^{\frac{2}{1-\frac{\nu}{2}}}\Big/[\epsilon^{\frac{1}{1-\frac{\nu}{2}}}(\ln\frac{1}{\rho})^{\frac{2}{1-\frac{\nu}{2}}}]\right)$.
\end{theorem}
When $\epsilon\ll\ln\frac{1}{\rho}$ and $0<\nu<1$, the result in Theorem \ref{th1a} significantly improves the speed of TD.  When $\nu=1$, the complexity is the same as \eqref{complexity-ori}.

We explain a little about assumption \eqref{fastass}. Note that the boundedness of the sequence $(\ttheta^k)_{k\geq 0}$ together with Assumptions 1 and 3 directly yields $v^k=O(k)$ (i.e., $\nu=1$).
In fact, assumption \eqref{fastass} is very standard for analyzing adaptive stochastic optimization   \cite{liao2021local,duchi2011adaptive,reddi2019convergence,chen2018universal,chen2018convergence,liu2019towards}.   As far as we know, there is no very good explanation of the superiority of the adaptive SGD  rather than assumption \eqref{fastass}; it is widely used in previous works because many training tasks enjoy sparse stochastic gradients.

In our AdaTD algorithm, we have
$$ {\bf g}^k=[r(s_{k+1},s_k)+\gamma V_{\ttheta^k}(s_{k+1})-V_{\ttheta^k}(s_k)]\phi(s_{k}),$$
which keeps the sparsity of $\phi(s_{k})$ because $r(s_{k+1},s_k)+\gamma V_{\ttheta^k}(s_{k+1})-V_{\ttheta^k}(s_k)\in\RR$.
Then if the features are  sparse, the stochastic semi-gradients are also  sparse.
In RL, a class of methods, called  as state-aggregation based approaches  \cite{mendelssohn1982iterative,bertsekas1989adaptive,abel2016near,duan2019state,li2006towards,van2006performance}, usually uses  sparse features to approximate Bellman equation linearly.
The main idea of  state-aggregation is  dividing the whole state apace into a few mutually disjoint   clusters (i.e., $\mathcal{S}=\bigcup_{i}^{\tilde{d}} \hat{s}_i$ and $\hat{s}_i\bigcap \hat{s}_j=\emptyset$ if $i\neq j$) and regards each cluster as a meta-state.
Compared with the vanilla MDP with linear approximation,
 the state-aggregation based one  employs a structured feature matrix
 $$\tilde{\Phi}:=[\phi(s_1), \phi(s_2),\ldots, \phi(s_{|\mathcal{S}|})]^{\top}\in\RR^{{|\mathcal{S}|} \times \tilde{d}},$$
where $\phi(s_i)[j]=1$ if $s_i\in \hat{s}_j$ and $\phi(s_i)[j]=0$ if $s_i\notin \hat{s}_j$. We can see that $\phi(s_i)$ has only one non-zero element, which is very sparse. And then, the stochastic semi-gradients in the AdaTD used to the state-aggregation MDP are very sparse.

State aggregation is a degenerated form of  linear representations. More general is the tile coding, which still retains the sparsity structure
\cite{sutton1996generalization}.

\section{Extension to   Projected  Adaptive TD($\lambda$)}
This section contains the adaptive TD($\lambda$) algorithm and its finite-time convergence analysis.

\subsection{Algorithm development}
Using existing analysis of TD($\lambda$) \cite{sutton1998introduction,tsitsiklis1996analysis}, if $V_\mu $ solves the Bellman equation \eqref{eq.bell}, it also solves
$$\Tc_\mu^{(m)} V_\mu = V_\mu,~m\in \mathbb{Z}^+,$$ where $\Tc_\mu^{(m)}$ denotes the $m$-step of $\Tc_\mu$. In this case, we can also represent $V_\mu$  as
$$V_\mu(s) = \mathbb{E}\big[\sum_{t=0}^m \gamma^t \Rc(s_t)+\gamma^{m+1}V_{\mu}(s_{m+1}) \mid s_0=s \big].$$
Given $\lambda\in [0,1]$ and
$V_\mu(s)=(1-\lambda)\sum_{m=0}^{\infty}\lambda^{m}V_\mu(s),$
 $\forall  s\in \Sc$, the $\lambda$-averaged Bellman operator is given by
\begin{align}\label{eq.bell-oprt2}
&(\Tc_\mu^{\lambda} V)(s)= (1-\lambda)\sum_{m=0}^{\infty}\lambda^{m}\nonumber\\
&\quad\times\mathbb{E}\left[ \sum_{t=0}^m \gamma^t \Rc(s_t)+\gamma^{m+1}V(s_{m+1}) \biggiven s_0=s \right].
\end{align}
Comparing operator $\Tc_\mu$ and the $\lambda$-averaged Bellman operator, it is clear that $\Tc_\mu^{0}=\Tc_\mu$.

By defining
$$
	{\bf g}^{\lambda}(\ttheta)=\Phi^{\top}\textrm{Diag}(\pi)(T_{\mu}^{\lambda}\Phi\ttheta-\Phi\ttheta),
$$
 the stochastic update in TD($\lambda$) can be presented as \eqref{sampling2}.
Similar to TD(0), it has been established in \cite{tsitsiklis1996analysis} and \cite{van1998learning} that
$$
	\lim_{k\rightarrow\infty}\EE[ \overline{{\bf g}}^\lambda(\ttheta;s_{k},s_{k+1},{\bf z}^k)]={\bf g}^{\lambda}(\ttheta).
$$
Like the limiting update direction ${\bf g}(\ttheta)$ in TD(0), a critical property of the update direction in TD($\lambda$) is given by
\begin{align}
\langle\ttheta^*-\ttheta, {\bf g}^{\lambda}(\ttheta)\rangle\geq(1-\alpha)\omega\|\ttheta-\ttheta^*\|^2,
\end{align}
where $\alpha:= \gamma(1-\lambda)/(1-\gamma\lambda)$ for any $\ttheta\in\RR^d$.
By denoting ${\bf z}^{\infty}:=\sum_{t=0}^{\infty}(\gamma\lambda)^t \hat{s}_{t}$, where $(\hat{s}_{1},\hat{s}_{2},\ldots)$ is the stationary sequence. Then, it also holds
$$
	\EE [\overline{{\bf g}}^\lambda(\ttheta;\hat{s}_{k},\hat{s}_{k+1},{\bf z}^{\infty})]={\bf g}^{\lambda}(\ttheta).
$$
We present AdaTD($\lambda$) in Algorithm \ref{alg2}. AdaTD($\lambda$) and AdaTD(0) differ in a different update protocol.
We directly have the following bounds for  AdaTD($\lambda$)
\begin{align*}\label{bound2}
\|\ttheta^{k}-\ttheta^{*}\|\leq 2\hat{R},~~~\|{\bf g}^k\|  \leq \hat{G},~~~\|{\bf m}^k\|  \leq \hat{G}
\end{align*}
where  $\hat{G}:=2\hat{R}+B$.

\begin{figure*}[t]
\centering
\begin{tabular}{cccc}
\hspace{-0.55cm}
\includegraphics[width=0.25\textwidth]{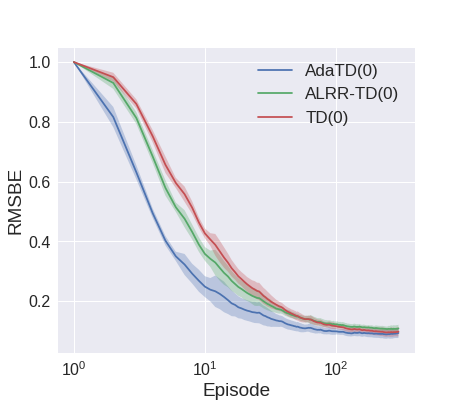}
\hspace{-0.5cm}
\includegraphics[width=0.25\textwidth]{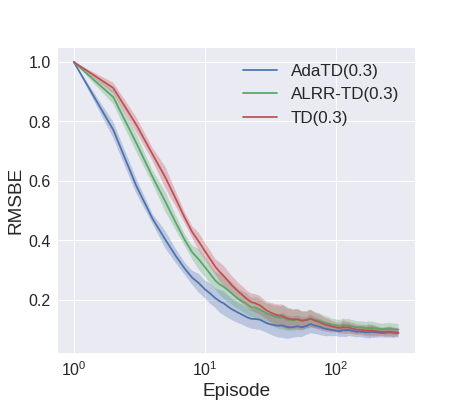}
\hspace{-0.5cm}
\includegraphics[width=0.25\textwidth]{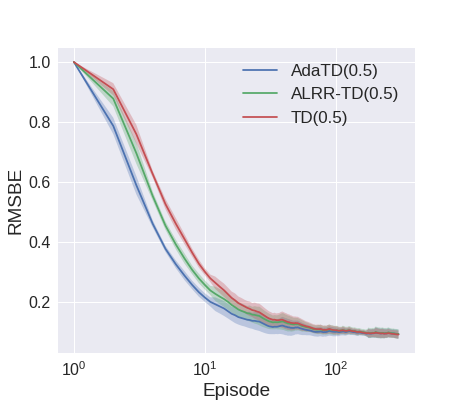}
\hspace{-0.5cm}
\includegraphics[width=0.25\textwidth]{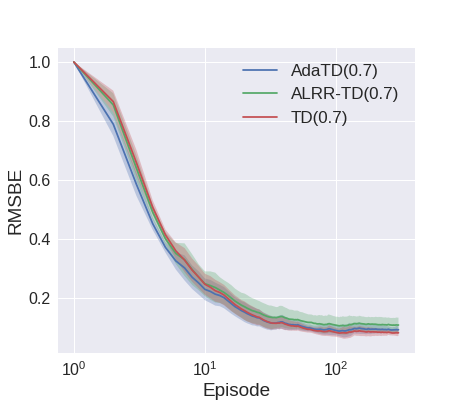}
\end{tabular}
\vspace*{-0.2cm}
  \caption{Runtime mean squared Bellman error (RMSBE) versus episode in `Mountain Car' when $\lambda = 0$, $0.3$, $0.5$ and $0.7$.}
\label{fig:rltest1}
\vspace*{-0.3cm}
\end{figure*}
\begin{figure*}[t]
\centering
\begin{tabular}{cccc}
\hspace{-0.55cm}
\includegraphics[width=0.25\textwidth]{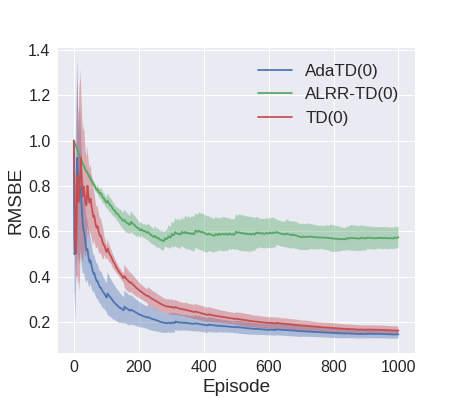}
\hspace{-0.5cm}
\includegraphics[width=0.25\textwidth]{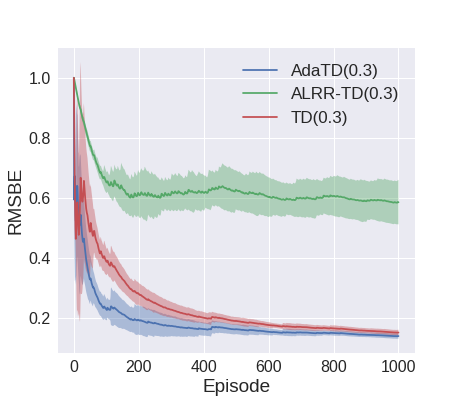}
\hspace{-0.5cm}
\includegraphics[width=0.25\textwidth]{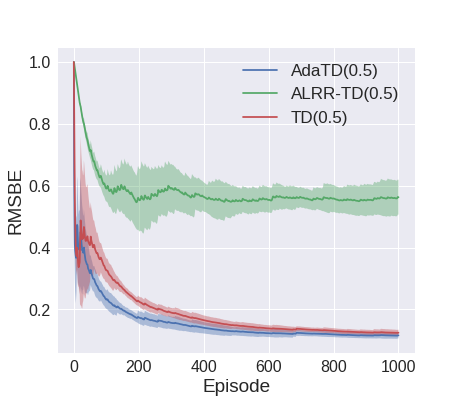}
\hspace{-0.5cm}
\includegraphics[width=0.25\textwidth]{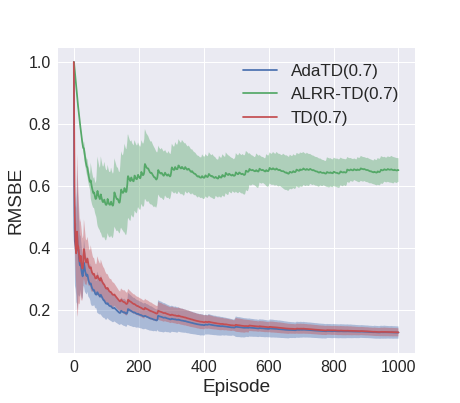}
\end{tabular}
\vspace*{-0.2cm}
  \caption{RMSBE versus episode in `Acrobot' when $\lambda = 0$, $0.3$, $0.5$ and $0.7$.}
\label{fig:rltest2}
\vspace*{-0.3cm}
\end{figure*}

\begin{figure*}[t]
\centering
\begin{tabular}{cccc}
\hspace{-0.55cm}
\includegraphics[width=0.25\textwidth]{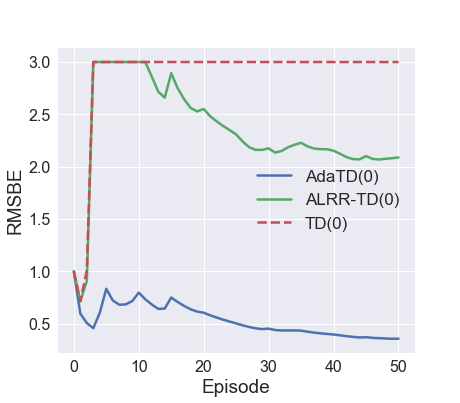}
\hspace{-0.5cm}
\includegraphics[width=0.25\textwidth]{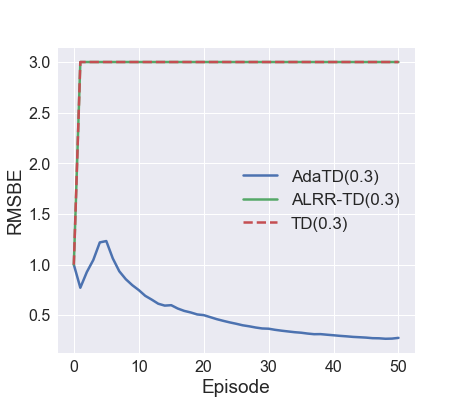}
\hspace{-0.5cm}
\includegraphics[width=0.25\textwidth]{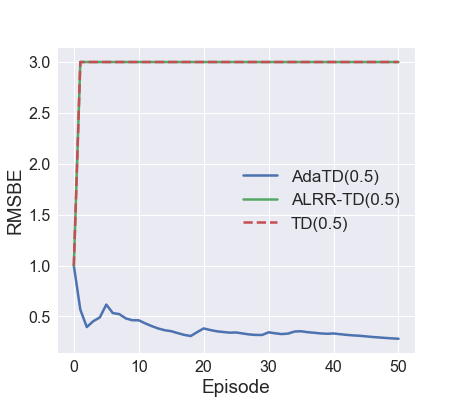}
\hspace{-0.5cm}
\includegraphics[width=0.25\textwidth]{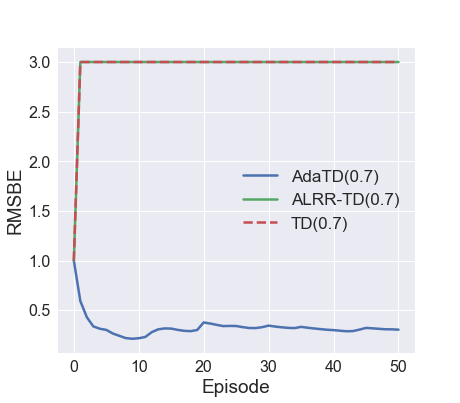}
\end{tabular}
\vspace*{-0.2cm}
  \caption{RMSBE versus episode in `Acrobot' when step size is large. The maximum RMSBE is capped since both ALRR and vanilla TD diverge.}
\label{fig:rltest5}
\vspace*{-0.3cm}
\end{figure*}

 \subsection{Finite-time analysis of projected AdaTD($\lambda$)}
The analysis of TD($\lambda$) is more complicated than   TD due to the existence of ${\bf z}^k$.
 To this end, we need to bound the sequence $(\EE{\bf z}^k)_{k\geq 0}$ in Lemma \ref{boundzk}.
 Similar to the analysis of  AdaTD(0) in Sec. 3, we need to consider the  ``delayed" expectation.
  For a fixed $K_0$, we consider the following error decomposition
  \begin{align*}
       & \EE [ \overline{{\bf g}}^\lambda(\ttheta^{k-K_0};s_{k},s_{k+1},{\bf z}^k)]=\EE [ \overline{{\bf g}}^\lambda(\ttheta^{k-K_0};s_{k},s_{k+1},{\bf z}^{\infty})]\\
  &+\EE  [\overline{{\bf g}}^\lambda(\ttheta^{k-K_0};s_{k},s_{k+1},{\bf z}^k)]-\EE  [\overline{{\bf g}}^\lambda(\ttheta^{k-K_0};s_{k},s_{k+1},{\bf z}^{\infty})].
  \end{align*}
 Therefore, our problem becomes bounding the difference between $\EE  [\overline{{\bf g}}^\lambda(\ttheta^{k-K_0};s_{k},s_{k+1},{\bf z}^{\infty})]$ and ${\bf g}^{\lambda}$, where the proof is similar to  Lemma \ref{legeo}. We can also establish the Lemma \ref{legeo2}.
We present the convergence of AdaTD($\lambda$) as follows.
\begin{theorem}\label{th2}
Suppose $(\ttheta^k)_{k\geq 0}$ are generated by AdaTD($\lambda$) with
$$
	\hat{R}\geq2B/(\sqrt{\omega}(1-\gamma)\sqrt{1-\frac{\gamma(1-\lambda)}{1-\gamma\lambda}})
$$
under the Markovian observation. Given  $\eta>0, \delta>0, 0\leq \beta<1$, $0\leq\lambda\leq 1$,  we have
\begin{align*}
\min_{1\leq k\leq K}\EE(\|\ttheta^*-\ttheta^k\|^2)&\leq \Big(C_1^{\lambda}\ln(\frac{\delta+KG^2}{\delta})\Big)/\sqrt{K}\\
&+C_2^{\lambda}/\sqrt{K},
\end{align*}
where $C_1^{\lambda}$ and  $C_2^{\lambda}$  are   given as
\begin{align*}
 &C_1^{\lambda}:=\frac{16(\ln K/\ln\frac{1}{\rho})^2\hat{G}}{\delta^{1/2}(1-\gamma)(1-\alpha)\omega^2}+\frac{2\eta\beta\hat{G}}{(1-\beta)(1-\alpha)\omega}\\
 &\quad+\frac{\eta\hat{G}}{(1-\alpha)\omega}+\frac{4R\hat{G}^2}{\delta(1-\alpha)\omega}=O(\frac{(\ln K/\ln\frac{1}{\rho})^2}{\sqrt{\delta}(1-\gamma\lambda)}),\\
 &C_2^{\lambda}:=\frac{4\hat{R}^2\hat{G}}{\eta(1-\alpha)\omega}+\frac{4\hat{R}\hat{G}^2\eta}{\delta^{\frac{1}{2}}(1-\beta)(1-\alpha)\omega}\\
 &\quad+\frac{2\hat{R}\hat{G}^2\sum_{k=1}^{+\infty}\zeta_{k}}{\sqrt{\delta}(1-\gamma\lambda)(1-\alpha)\omega}+\frac{4\hat{R} \bar{\kappa}(B+\hat{R}\gamma+\hat{R})\hat{G} }{(1-\alpha)\sqrt{\delta}\omega}\\
 &\quad=O(\frac{1}{\sqrt{\delta}(1-\beta)}+\frac{\max\{\gamma\lambda,\rho\}  }{\sqrt{\delta}(1-\gamma\lambda)(1-\max\{\gamma\lambda,\rho\})^2}),
\end{align*}
and  $\zeta_k:=\sum_{t=1}^{k}(\gamma\lambda)^{k-t}\rho^t$.
\end{theorem}
If $\lambda=0$, it then holds $\zeta_k\equiv 0$.\footnote{For convenience, we follow the convention $0^0=0$.}
 It is also easy to see $\zeta_k\leq  k(\max\{\gamma\lambda,\rho\})^{k}$.
But we do not want to use $k(\max\{\gamma\lambda,\rho\})^{k}$ to replace $\zeta_k$ in Theorem \ref{th2} because the bound will not diminish when $\lambda=0$.
When $\lambda=0$, Theorem \ref{th2} reduces to Theorem \ref{th1}. Using similar argument like Theorem \ref{th1}, to achieve $\epsilon$ error for $\min_{1\leq k\leq K} \EE(\|\ttheta^*-\ttheta^k\|^2)$, the numer of iteration $K$ is again $\tilde{O}( \ln^4\frac{1}{\epsilon}/(\epsilon^2\ln^4\frac{1}{\rho}))$.

In  \cite{van1998learning} and \cite{bertsekas2011dynamic}, it has been proved that the zero point of ${\bf g}^{\lambda}(\ttheta)$ (i.e., the limiting point of  projected AdaTD($\lambda$)) satisfies
\begin{equation}\label{dislam}
\|\Phi \ttheta^*-V_{\mu}\|_{\pi}\leq \frac{1}{\sqrt{1-[\frac{\gamma(1-\lambda)}{1-\gamma\lambda}]^2}}\|\PP_{\Phi}(V_{\mu})-V_{\mu}\|_{\pi},
\end{equation}
where $V_{\mu}$ is the solution of the Bellman equation, and $\PP_{\Phi}$ is the projection  onto the span of $\Phi$'s columns, and $\|{\bf V}\|_{\pi}:=\sqrt{\sum_{i=1}^{|\mathcal{S}|}\pi(s_i)[V(i)]^2}$. From \eqref{dislam}, we can see that   the limiting point of projected AdaTD($\lambda$) yields  a better approximator  than  that of  projected AdaTD($0$) when $0<\lambda\leq 1$;  additionally,  $\lambda=1$  provides the best approximation gap.
Nevertheless, only the approximation performance is insufficient   to force us to use $\lambda=1$ in all scenarios because $\lambda$ also affects  the  convergence speed.
We  explain why non-zero $\lambda$ hurts the iterative speed of  projected  AdaTD($\lambda$):
When $\lambda\neq 0$, $C_2^{\lambda}$ has an extra term than $C_2$, i.e.,
$\frac{2\hat{R}\hat{G}^2\sum_{k=1}^{+\infty}\zeta_{k}}{\sqrt{\delta}(1-\gamma\lambda)(1-\alpha)\omega}=\Theta(\frac{\max\{\gamma\lambda,\rho\}  }{(1-\gamma\lambda)(1-\max\{\gamma\lambda,\rho\})^2})$,
which increases when $\lambda$ increases.
In summary, $\lambda$ reflects the trade-off between accuracy of the limiting point and convergence speed.

\medskip

The acceleration of AdaTD($\lambda$) is also provable when the semi-gradients are ``sparse".
\begin{proposition}\label{pro1}
Suppose the conditions of Theorem \ref{th2} hold. To achieve $\epsilon$-accuracy for $\min_{1\leq k\leq K}\{\EE(\|\ttheta^*-\ttheta^k\|^2)\} $,  the needed iteration is
$\tilde{O}( \ln^4\frac{1}{\epsilon}/(\epsilon^2\ln^4\frac{1}{\rho}))$. Further if
\eqref{fastass}
holds, the complexity can be improved as
 \begin{equation}
    \tilde{O}\left((\ln\frac{1}{\epsilon})^{\frac{2}{1-\frac{\nu}{2}}}\Big/[\epsilon^{\frac{1}{1-\frac{\nu}{2}}}(\ln\frac{1}{\rho})^{\frac{2}{1-\frac{\nu}{2}}}]\right).
 \end{equation}
\end{proposition}

\begin{figure*}[t]
\centering
\begin{tabular}{cccc}
\hspace{-0.55cm}
\includegraphics[width=0.25\textwidth]{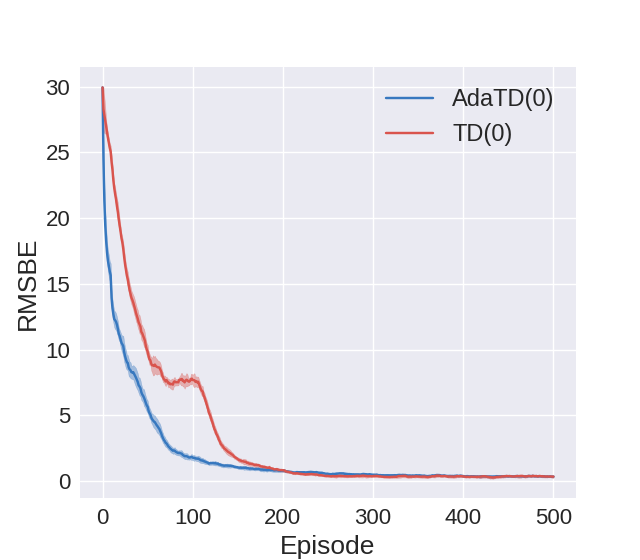}
\hspace{-0.5cm}
\includegraphics[width=0.25\textwidth]{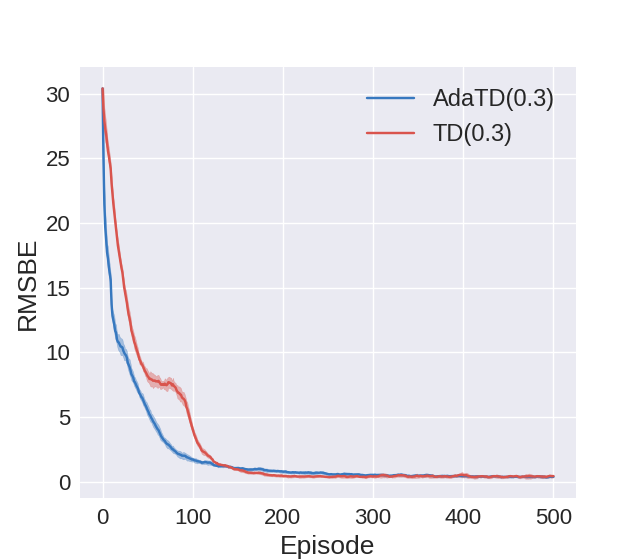}
\hspace{-0.5cm}
\includegraphics[width=0.25\textwidth]{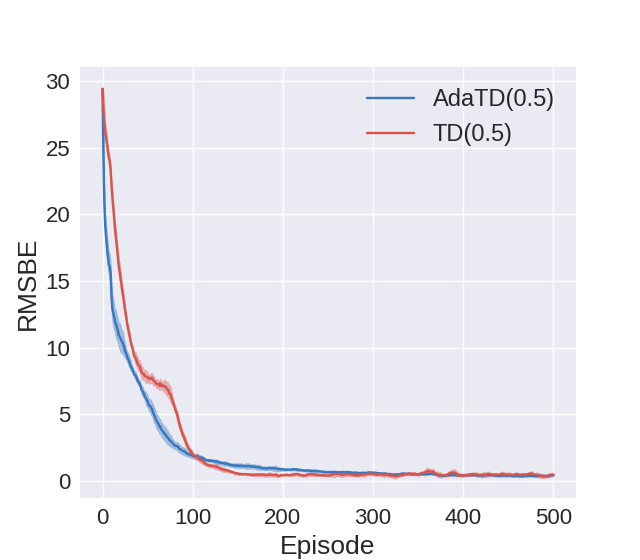}
\hspace{-0.5cm}
\includegraphics[width=0.25\textwidth]{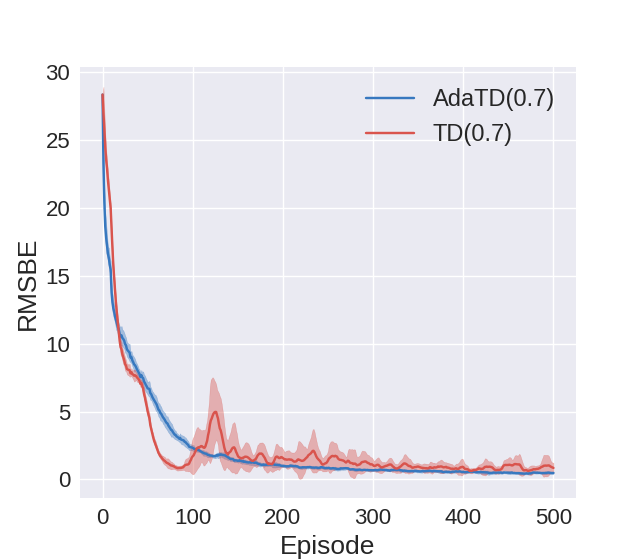}
\end{tabular}
\vspace*{-0.2cm}
  \caption{RMSBE versus episode in `CartPole' when $\lambda = 0$, $0.3$, $0.5$ and $0.7$.}
\label{fig:rltest3}
\vspace*{-0.2cm}
\end{figure*}
\begin{figure*}[t]
\centering
\begin{tabular}{cccc}
\hspace{-0.55cm}
\includegraphics[width=0.25\textwidth]{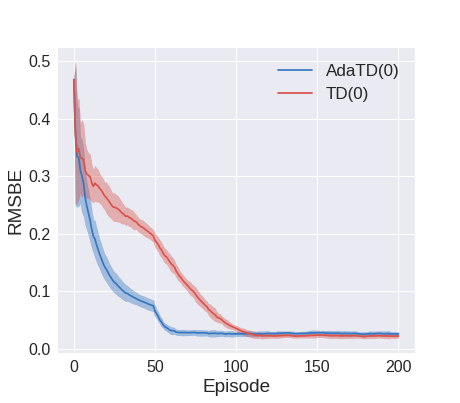}
\hspace{-0.5cm}
\includegraphics[width=0.25\textwidth]{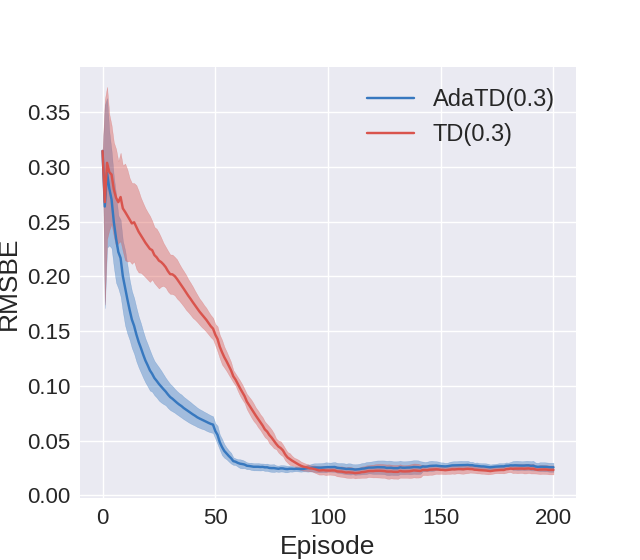}
\hspace{-0.5cm}
\includegraphics[width=0.25\textwidth]{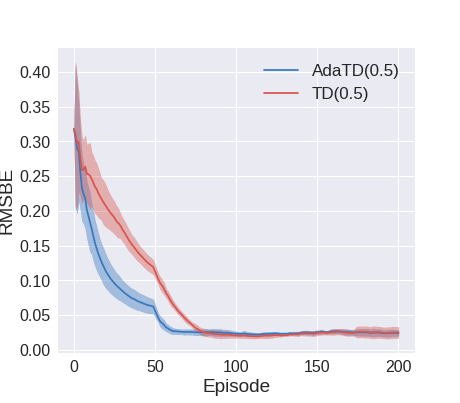}
\hspace{-0.5cm}
\includegraphics[width=0.25\textwidth]{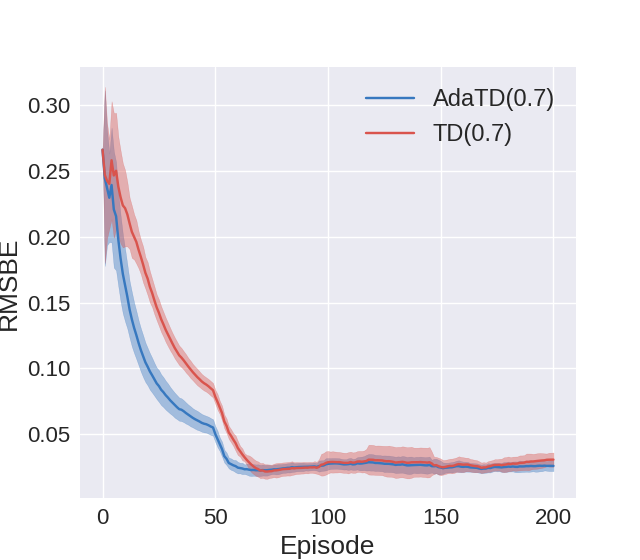}
\end{tabular}
\vspace*{-0.2cm}
  \caption{RMSBE versus episode in `Navigation' when $\lambda = 0$, $0.3$, $0.5$ and $0.7$.}
\label{fig:rltest4}
\vspace*{-0.3cm}
\end{figure*}

\begin{figure*}[t]
\centering
\begin{tabular}{cccc}
\hspace{-0.55cm}
\includegraphics[width=0.25\textwidth]{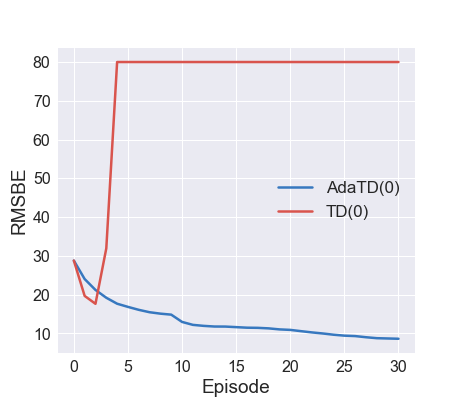}
\hspace{-0.5cm}
\includegraphics[width=0.25\textwidth]{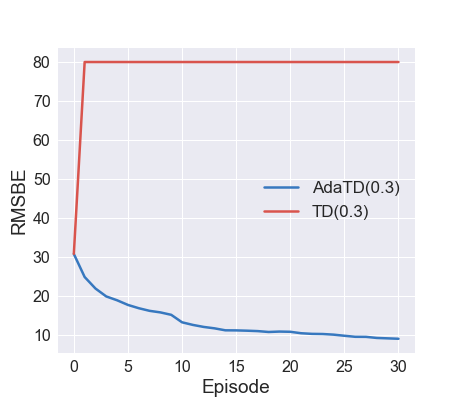}
\hspace{-0.5cm}
\includegraphics[width=0.25\textwidth]{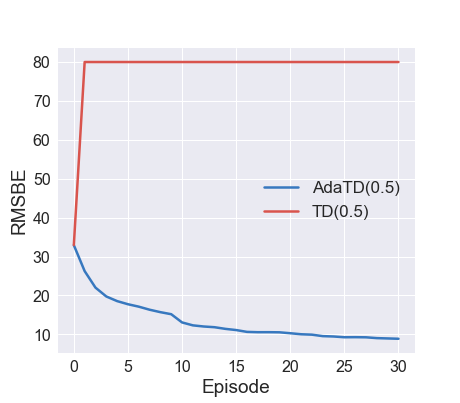}
\hspace{-0.5cm}
\includegraphics[width=0.25\textwidth]{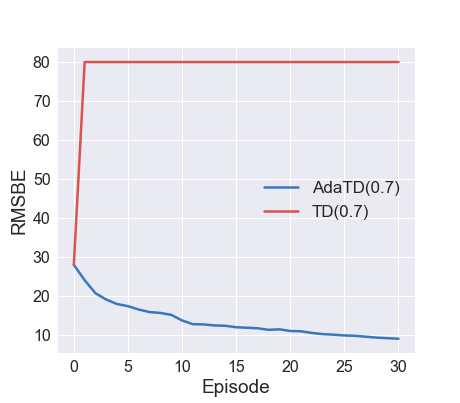}
\end{tabular}
\vspace*{-0.2cm}
  \caption{RMSBE versus episode in `Cartpole' when step size is large. The maximum RMSBE is capped.}
\label{fig:rltest6}
\vspace*{-0.3cm}
\end{figure*}


\begin{figure*}[t]
\centering
\begin{tabular}{cccc}
\hspace{-0.55cm}
\includegraphics[width=0.25\textwidth]{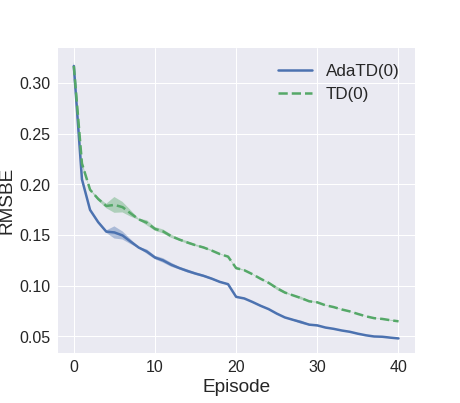}
\hspace{-0.5cm}
\includegraphics[width=0.25\textwidth]{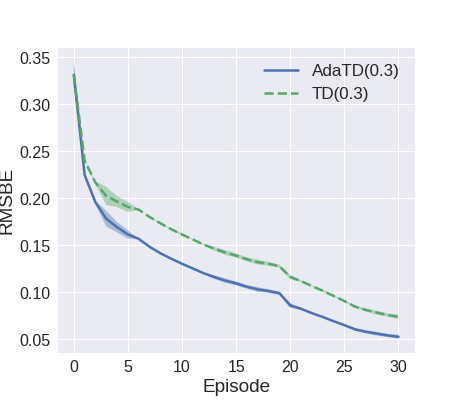}
\hspace{-0.5cm}
\includegraphics[width=0.25\textwidth]{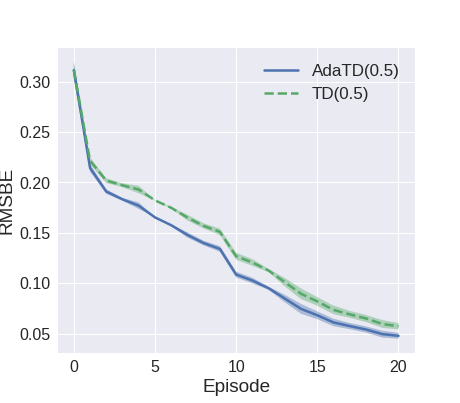}
\hspace{-0.5cm}
\includegraphics[width=0.25\textwidth]{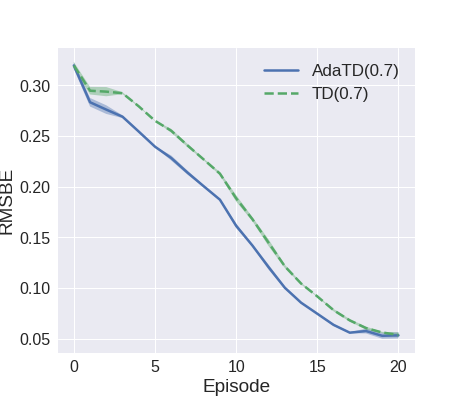}
\end{tabular}
\vspace*{-0.2cm}
  \caption{RMSBE versus episode in random MDP with sparse state features.}
\label{fig:rltest7}
\vspace*{-0.3cm}
\end{figure*}


\section{Numerical Simulations}

To validate the analysis and show the effectiveness of our algorithms, we tested AdaTD(0) and AdaTD($\lambda$) on several commonly used RL tasks.
As briefly highlighted below, the first three tasks are from the popular OpenAI Gym \cite{openai2016}, and the fourth task is a single-agent version of the \emph{Navigation} task in \cite{lowe2017nips}.

We compared our algorithms with other policy evaluation methods using the runtime mean squared Bellman error (RMSBE). In each test,  the policy is the same for all the algorithms when the value parameter is updated separately.
In the first two tasks, the value function is approximated using linear functions.
In the last two tasks, the value function is parameterized by a neural network.
In the linear tasks, for different values of $\lambda$, we compared AdaTD($\lambda$) algorithm, the TD($\lambda$), and ALRR algorithm in \cite{gupta2019finite}.
For fair comparison, we changed the update step in the original ALRR algorithm to a single time scale TD($\lambda$) update.
In the non-linear tasks, we extended our algorithm to non-linear cases and compared it with TD($\lambda$).
Since ALRR was not designed for the neural network-parameterized cases, we did not include it in the non-linear TD tests. In all tests, the curves are generated by taking the average of 10 Monte-Carlo runs.

\subsection{Experiment Settings}

\noindent\textbf{Mountain Car.}
Algorithms were tested when $\lambda$ = $0$, $0.3$, $0.5$ and $0.7$.
For all methods, we set max episode = 300, batch size (horizon) = 16. In vanilla TD method, $\eta$ = $0.7$.
In ALRR-TD($\lambda$), $\eta_0$ = $1.0$, $\sigma$ = $0.001$, $\xi$ = $1.2$.
In AdaTD($\lambda$), $\eta_0$ = $3.0$, $\delta$ = $0.01$, $\beta$ = $0.3$.

\noindent\textbf{Acrobot.}
In all three methods, max episode = 1000 and batch size = 48.
In vanilla TD($\lambda$), $\eta$ = $0.05$ when $\lambda = 0$, otherwise $\eta$ = $0.04$.
In ALRR-TD($\lambda$), $\eta_0$ = $0.001$, $\sigma$ = $0.001$, $\xi$ = $1.2$.
In AdaTD($\lambda$), $\eta_0$ = $9$, $\delta$ = $1$. When $\lambda = 0$ and $0.3$, $\beta$ = $0.5$, otherwise $\beta$ = $0.9$.

\noindent\textbf{CartPole.}
We used a neural network to approximate the value function. The neural network has two hidden layers each with 128 neurons and ReLU activation.
For both methods, we set max episode = 500 and batch size = 32.
In TD($\lambda$), $\eta$ = $0.02$.
In AdaTD($\lambda$), $\eta_0$ = $1.5$, $\delta$ = $0.01$, $\beta$ = $0.2$.

\noindent\textbf{Navigation.}
We used a neural network to approximate the value function. The neural network has two hidden layers each with 64 neurons and ReLU activation.
For both methods, we set max episode = 200 and batch size = 20.
In TD($\lambda$), $\eta = 0.2$.
In AdaTD($\lambda$), $\eta_0$ = $0.7$, $\delta$ = $0.01$, $\beta$ = $0.2$.

\subsection{Numerical Results}
In the test of Mountain Car, the performance of all three methods is close, while AdaTD($\lambda$) still has a small advantage over other two when $\lambda$ is small.
In the Acrobot task, the initial step size is relatively large for AdaTD($\lambda$), but AdaTD($\lambda$) is able to adjust the large initial step size and guarantee afterwards convergence.
Note there is a major fluctuation in average loss around episode $30$. TD($\lambda$) has constant step size, and thus it is more vulnerable to the fluctuation than AdaTD($\lambda$). As a result, our algorithm demonstrates better overall convergence speed over TD($\lambda$).

In Figures \ref{fig:rltest3} and \ref{fig:rltest4}, we tested our algorithm with neural network parameterization. In these tests,
the step size of TD($\lambda$) cannot be large due to stability issues. As a result, TD($\lambda$) is outperformed by AdaTD($\lambda$) where large step size is allowed. In fact, when $\lambda$ is large, a small step size of $\eta = 0.02$ still cannot guarantee the stability of TD($\lambda$). It can be observed in Figure \ref{fig:rltest3} that when $\lambda$ gets larger, i.e. the gradient magnitude is larger, original TD($\lambda$) becomes less stable. In comparison, AdaTD($\lambda$) has exhibited robustness to the choice of $\lambda$ and the large initial step size in this test.

We conduct two experiments in Figures \ref{fig:rltest5} and \ref{fig:rltest6} where we deliberately skip the step size tuning process and applies large step size to all methods. We show that our algorithm is more robust to large step size and therefore is easier to tune in practice.
The hyper-parameters of the two tests are the same as before except for step sizes. In Acrobot, step size is set to 1. In Cartpole, step size is set to 0.8.

In Figure
\ref{fig:rltest7}, we have also provided another test to evaluate the performance of AdaTD under sparse state features. The MDP has a discrete state space of size 50, and a discrete action space of size 4. The transition matrix and reward table are randomly generated with each element in $(0,1)$. To ensure the sparsity of gradients, we construct sparse state features as described in Section \ref{section:analysis_proj_atd}.
In this test, we select $\eta=0.45$ for TD, and $\eta_0=0.5$,$\delta$ = $1.0$, $\beta$ = $0.5$ for AdaTD.
\section{Proofs}

\subsection{Proofs for AdaTD(0)}

This part contains the proofs of the main results and leaves the proofs of the technical lemmas in the appendix.
\subsubsection{Technical Lemmas}
\begin{lemma}[\cite{chen2018convergence,li2018convergence}]\label{le1}
For $0\leq a_i\leq \bar{a}$ and $\delta>0$, we have
\begin{align*}
\sum_{t=1}^{K_0} \frac{a_t}{\delta+\sum_{i=1}^t a_i}\leq\ln(\delta+K_0\bar{a})-\ln\delta.
\end{align*}
\end{lemma}

\begin{lemma}[\cite{tsitsiklis1996analysis}]\label{boundg}
For any $\ttheta\in\RR^d$, it follows
\begin{align}\label{bg}
\langle \ttheta^*-\ttheta, {\bf g}(\ttheta)\rangle\geq (1-\gamma)\omega\|\ttheta^*-\ttheta\|^2.
\end{align}
\end{lemma}

In the proofs, we use three shorthand notation for simplifications.
Those three notation are all related to the iteration $k$. Assume  $({\bf m}^k)_{k\geq 0}$, $(\ttheta^k)_{k\geq 0}$,  $(v^k)_{k\geq 0}$ are all generated by AdaTD(0). We denote
\begin{align}\label{notation1}
   \begin{cases}
    \hat{\bf m}_k&:=  \EE \Big[\|{\bf m}^{k}\|^2/(v^{k}+\delta)\Big],\\
    \Delta_k&:=\EE\left(\langle \ttheta^{k}-\ttheta^{*},  {\bf m}^k\rangle/(  v^k+\delta)^{\frac{1}{2}}\right),\\
   \phi_k&:=\EE\left(\langle \ttheta^{*}-\ttheta^{k},  {\bf g}(\ttheta^k)\rangle/(  v^{k-1}+\delta)^{\frac{1}{2}}\right)\\
\Re_k&:= \frac{8K_0(1-\beta)}{\delta^{1/2}(1-\gamma)\omega}\sum_{h=K_0}^1\hat{\bf m}_{k-h}+\eta\beta\hat{\bf m}_{k}\\
&\quad+2RG(1-\beta)\frac{1}{\delta} \cdot\frac{\|{\bf g}^k\|^2}{v^k+\delta}+\frac{2R(1-\beta)\kappa \rho^{K_0}}{\sqrt{\delta}}.
   \end{cases}
\end{align}
The above notations will be used in the remaining lemmas.

\begin{lemma}\label{core1}
Let $(\Delta_k)_{k\geq 0}$  and $(\Re_k)_{k\geq 0}$ be defined in \eqref{notation1}, the following result holds for AdaTD(0)
\begin{align}
\Delta_k+\frac{(1-\beta)}{2}\phi_k\leq\beta\Delta_{k-1}+\Re_k.
\end{align}
On the other hand, we can bound $\Delta_k$ as
$|\Delta_k|\leq \eta\frac{2RG}{\delta^{\frac{1}{2}}}.$
\end{lemma}


\subsubsection{Proof of Theorem \ref{th1}}
Given $K,K_0\in \mathbb{Z}^+$ and $K\geq K_0-1$, Lemma \ref{core1} tells
\begin{align*}
\Delta_k+\frac{1-\beta}{2}\phi_k\leq \beta\Delta_{k-1}+\Re_k,
\end{align*}
where $K_0+1\leq k\leq K$.
Summing together, we then get
\begin{align}\label{th1-t1}
&\frac{(1-\beta)}{2}\sum_{k=K_0+1}^K \phi_k \nonumber\\ &\leq-\Delta_K+(\beta-1)\sum_{k=K_0}^{K-1}\Delta_k+\sum_{k=K_0+1}^K\Re_k\nonumber\\
&\leq (\beta-1)\sum_{k=K_0}^{K-1}\Delta_k+\sum_{k=K_0+1}^K\Re_k+\eta\frac{2RG}{\delta^{\frac{1}{2}}}.
\end{align}
With direct calculations, we get
\begin{align*}
&\|\ttheta^*-\ttheta^{k+1}\|^2\\
&\leq\|\ttheta^*-\ttheta^{k}-\eta {\bf m}^k/\sqrt{v^k+\delta}\|^2\\
&\leq \|\ttheta^{*}-\ttheta^{k}\|^2+\frac{2\eta \langle{\bf m}^k,\ttheta^{k}-\ttheta^{*}\rangle}{(  v^k+\delta)^{\frac{1}{2}}}+\frac{\eta^2\|{\bf m}^k\|^2}{v^k+\delta}.
\end{align*}
Taking total condition expectation gives us
\begin{align*}
\EE \|\ttheta^*-\ttheta^{k+1}\|^2\leq \EE \|\ttheta^*-\ttheta^{k}\|^2+2\eta\Delta_k+\eta^2\hat{\bf m}_k.
\end{align*}
That is also
\begin{align*}
\sum_{k=K_0}^{K-1} -\Delta_k \leq   \frac{\EE\|\ttheta^*-\ttheta^{K_0}\|^2}{2\eta}+\frac{\eta}{2}\sum_{k=K_0}^{K-1}\hat{\bf m}_k.
\end{align*}
Combining \eqref{th1-t1}, we are then led to
\begin{align}\label{lat}
&\sum_{k=K_0+1}^K\phi_k\nonumber\\
& \leq 2\sum_{k=K_0}^{K-1}(-\Delta_k)+\frac{2}{1-\beta}\sum_{k=K_0+1}^K\Re_k+\frac{4RG\eta}{\delta^{\frac{1}{2}}(1-\beta)}\nonumber\\
&\leq \frac{\EE\|\ttheta^*-\ttheta^{K_0}\|^2}{\eta}+\eta\sum_{k=K_0}^{K-1}\hat{\bf m}_k\nonumber\\
&\qquad+\sum_{k=K_0+1}^K\frac{2\Re_k}{1-\beta}+\frac{4RG\eta}{\delta^{\frac{1}{2}}(1-\beta)}\nonumber\\
&\leq \frac{4R^2}{\eta}+\eta\sum_{k=K_0}^{K-1}\hat{\bf m}_k+\sum_{k=K_0+1}^K\frac{2\Re_k}{1-\beta}+\frac{4RG\eta}{\delta^{\frac{1}{2}}(1-\beta)}.
\end{align}
Now, we turn to bound the right side of \eqref{lat}: with the definition of $\Re_k$,
\begin{align*}
&\eta\sum_{k=K_0}^{K-1}\hat{\bf m}_k+\frac{2}{1-\beta}\sum_{k=K_0+1}^K\Re_k\nonumber\\
&\leq \eta\sum_{k=K_0}^{K-1}\hat{\bf m}_k+ \left(\frac{16K_0^2}{\delta^{1/2}(1-\gamma)\omega}+\frac{2\eta\beta}{1-\beta}\right)\sum_{k=K_0+1}^K\hat{\bf m}_{k} \nonumber\\
&\qquad+\frac{4RG}{\delta} \sum_{k=K_0+1}^K\frac{\|{\bf g}^k\|^2}{v^k+\delta}+\frac{4R \kappa \rho^{K_0}}{\sqrt{\delta}}(K-K_0). \nonumber\\
\end{align*}
With  Lemma \ref{core0}, the bound is further bounded by
\begin{align}\label{th1-t2}
& \left(\eta+\frac{16K_0^2}{\delta^{1/2}(1-\gamma)\omega}+\frac{2\eta\beta}{1-\beta}\right)\sum_{k=1}^K\hat{\bf m}_{k} \nonumber\\
&\qquad+\frac{4RG}{\delta} \sum_{k=1}^K\frac{\|{\bf g}^k\|^2}{v^k+\delta}+\frac{4R \kappa \rho^{K_0}}{\sqrt{\delta}}(K-K_0) \nonumber\\
&\overset{\textrm{Lemma}~\ref{core0}}{\leq} \left( \frac{16K_0^2}{\delta^{1/2}(1-\gamma)\omega}+\frac{2\eta\beta}{1-\beta}+\eta\right) \nonumber\\
&\qquad\times\ln(\frac{\delta+(K-1)G^2}{\delta}) \nonumber\\
 &\qquad+\frac{4RG}{\delta} \ln(\frac{\delta+KG^2}{\delta})+\frac{4R \kappa \rho^{K_0}}{\sqrt{\delta}}(K-K_0),
\end{align}
where we used the inequality $\sum_{k=1}^K\frac{\|{\bf g}^k\|^2}{v^k+\delta}\leq \ln(\frac{\delta+KG^2}{\delta})$ (Lemma \ref{le1}).
We set $K_0= \ulcorner\ln K/\ln\frac{1}{\rho}\urcorner$. It is easy to see that $K\geq\frac{9}{4}(K_0+1)$ as $K$ is large.
On the other hand, with Lemma \ref{boundg}, we can get
\begin{align}\label{th1-t3}
\sum_{k=K_0}^K\phi_k
&\geq \sum_{k=K_0}^K \frac{(1-\gamma)\omega\EE\|\ttheta^*-\ttheta^k\|^2}{[(k-1)G^2+\delta]^{\frac{1}{2}}}\nonumber\\
 &\geq 2\Big(\sqrt{(K-1)G^2+\delta}-\sqrt{(K_0-1)G^2+\delta}\Big)/G^2\nonumber\\
 &\qquad\!\times(1-\gamma)\omega\min_{K_0\leq k\leq K}\EE\|\ttheta^*-\ttheta^k\|^2\nonumber\\
 &\geq (1-\gamma)\omega\sqrt{K}/G\cdot\min_{K_0\leq k\leq K}\EE\|\ttheta^*-\ttheta^k\|^2,
\end{align}
where we used $K\geq \frac{9}{4}(K_0+1)$ to get
$$2(\sqrt{(K-1)G^2+\delta}-\sqrt{(K_0-1)G^2+\delta})\geq \sqrt{K}G.$$
In this case, we then derive
\begin{align*}
& (1-\gamma)\omega\min_{1\leq k\leq K}\EE\|\ttheta^*-\ttheta^k\|^2\\
&\leq(1-\gamma)\omega \min_{K_0\leq k\leq K}~\EE\|\ttheta^*-\ttheta^k\|^2\\
&\leq\frac{4R \kappa \rho^{K_0}G}{\sqrt{\delta}}\sqrt{K}+\Big( \frac{16K_0^2 G}{\delta^{1/2}(1-\gamma)\omega}+\frac{2G\eta\beta}{(1-\beta)}+\eta G\\
&+\frac{4RG^2}{\delta} \Big)\times\ln(\frac{\delta+KG^2}{\delta})\frac{1}{\sqrt{K}}+(\frac{4R^2G}{\eta}+\frac{4RG^2\eta}{\delta^{\frac{1}{2}}(1-\beta)})\frac{1}{\sqrt{K}}.
\end{align*}
By setting $K_0:=\ulcorner\ln K/\ln\frac{1}{\rho}\urcorner$ and defining the constants involved in the theorem as
given in Theorem \ref{th1}.
we then proved the results.

\subsubsection{Proof of Theorem \ref{th1a}}
The proof lies on a re-estimate of \eqref{th1-t3}. Under condition \eqref{fastass}, we can derive
\begin{align}\label{th1a-t3}
\sum_{k=K_0}^K\phi_k&\geq \sum_{k=K_0}^K \frac{(1-\gamma)\omega\EE\|\ttheta^*-\ttheta^k\|^2}{[v^{k-1}+\delta]^{\frac{1}{2}}}\nonumber\\
 &\geq \sum_{k=K_0}^K \frac{(1-\gamma)\omega}{[c (k-1)^{\nu}+\delta]^{\frac{1}{2}}}\!\times\min_{K_0\leq k\leq K}\EE\|\ttheta^*-\ttheta^k\|^2\nonumber\\
 &\geq \frac{K^{1-\frac{\nu}{2}} (1-\gamma)\omega}{2(1-\nu/2)\sqrt{c}}\cdot\min_{K_0\leq k\leq K}\EE\|\ttheta^*-\ttheta^k\|^2,
\end{align}
as $K\geq 2^{\frac{1}{1-\nu/2}}K_0$. Noticing that \eqref{th1-t2} still holds, which means
$$\min_{1\leq k\leq K}\EE\|\ttheta^*-\ttheta^k\|^2=O(\frac{K_0^2\ln K}{K^{1-\frac{\nu}{2}}}+\rho^{K_0}K^{\frac{\nu}{2}}).$$
To achieve the $O(\epsilon)$-accuracy for $\min_{1\leq k\leq K}\{\EE\|\ttheta^*-\ttheta^k\|^2\} $,  we need
\begin{align*}
\rho^{K_0}K^{\frac{\nu}{2}}\leq \epsilon,
       \frac{K_0^2\ln K}{K^{1-\frac{\nu}{2}}}\leq  \epsilon.
\end{align*}
We then have $K_0=\tilde{O}\left(\ln\frac{1}{\epsilon}/\ln\frac{1}{\rho}\right)$ and $K=\tilde{O}\left((\ln\frac{1}{\epsilon})^{\frac{2}{1-\frac{\nu}{2}}}\Big/[\epsilon^{\frac{1}{1-\frac{\nu}{2}}}(\ln\frac{1}{\rho})^{\frac{2}{1-\frac{\nu}{2}}}]\right)$.

\subsection{Proofs for AdaTD($\lambda$)}
The proof of AdaTD($\lambda$) is similar to AdaTD($0$) but with involved technical lemmas being modified. The complete proof is given in the appendix.

\section{Conclusions}
We developed an improved variant of the celebrated TD learning algorithm in this paper.
Motivated by the close link between TD(0) and stochastic gradient-based methods, we developed adaptive TD(0) and TD($\lambda$) algorithms.
The finite-time convergence analysis of the novel adaptive TD(0) and TD($\lambda$) algorithms has been established under the Markovian observation model.
While the worst-case convergence rates of Adaptive TD(0) and TD($\lambda$) are similar to those of the original TD(0) and TD($\lambda$), our numerical tests on several standard benchmark tasks demonstrate the effectiveness of our adaptive approaches. Future work includes variance-reduced and decentralized AdaTD approaches.

\newpage

\section{Proofs of Technical Lemmas for AdaTD(0)}
\subsection{Proof of Lemma \ref{legeo}}
Given a fixed integer $K_0$, with Assumptions \ref{ass0} and \ref{ass1},  by using a shorthand notation $$\Xi:=\left(\Rc(s,s')+\gamma \phi(s')^\top \ttheta^{k-K_0}- \phi(s)^\top \ttheta^{k-K_0} \right)\phi(s),$$ we have
\begin{align}\label{biasedexp}
&\EE( \overline{{\bf g}}(\ttheta^{k-K_0};s_{k},s_{k+1})\mid\sigma^{k-K_0})=\sum_{s, s' \in \Sc} \pi(s)\cdot\Pc(s'|s) \cdot \Xi \nonumber\\
& +\sum_{s, s' \in \Sc} (\Pc(s_k\mid s_{k-K_0}=s)-\pi(s)) \cdot\Pc(s'|s)\cdot\Xi.
\end{align}
Noticing that the following expectation
\begin{align*}
&\sum_{s, s' \in \Sc}  \pi(s)\cdot\Pc(s'|s) \cdot \Xi={\bf g}(\ttheta^{k-K_0}).
\end{align*}
With Assumption 2,   the second term in right side of \eqref{biasedexp} then can be bounded by $S\bar{\kappa}(B+R\gamma+R)\rho^{K_0}$. By setting $\kappa:=S\bar{\kappa}(B+R\gamma+R)$ and using $\EE(\EE(\cdot\mid\sigma))=\EE(\cdot)$, we then proved the result.

\subsection{Proof of Lemma \ref{lemmayu}}
With the fact
\begin{align*}
&{\bf g}^k= {\bf g}(\ttheta^{k})+{\bf g}^k-\overline{{\bf g}}(\ttheta^{k-K_0};s_{k},s_{k+1})\\ &+\overline{{\bf g}}(\ttheta^{k-K_0};s_{k},s_{k+1})-{\bf g}(\ttheta^{k-K_0}) +{\bf g}(\ttheta^{k-K_0})-{\bf g}(\ttheta^{k}),
\end{align*}
it then follows
\begin{align}\label{yu-le1}
&\EE\Big[\langle \ttheta^{k}-\ttheta^{*},  {\bf g}^k\rangle/( v^{k-1}+\delta)^{\frac{1}{2}}\Big]\nonumber\\
&=\EE\Big[\langle\ttheta^{k}-\ttheta^{*}, {\bf g}(\ttheta^{k})\rangle/(v^{k-1}+\delta)^{\frac{1}{2}}\Big]+\textrm{I}+\textrm{II}+\textrm{III},
\end{align}
where
\begin{align*}
     \begin{cases}
    \textrm{I}&:=\frac{\EE\mid\langle \ttheta^{k}-\ttheta^{*},  \left[ {\bf g}^k-\overline{{\bf g}}(\ttheta^{k-K_0};s_{k},s_{k+1})\right]\rangle\mid}{(v^{k-1}+\delta)^{\frac{1}{2}}},\\
    \textrm{II}&:=\frac{\EE\mid\langle \ttheta^{k}-\ttheta^{*},  \left[ \overline{{\bf g}}(\ttheta^{k-K_0};s_{k},s_{k+1})-{\bf g}(\ttheta^{k-K_0})\right]\rangle\mid}{(v^{k-1}+\delta)^{\frac{1}{2}}},\\
    \textrm{III}&:=\frac{\EE\mid\langle \ttheta^{k}-\ttheta^{*},  \left[ {\bf g}(\ttheta^{k-K_0})-{\bf g}(\ttheta^{k})\right]\rangle\mid}{(v^{k-1}+\delta)^{\frac{1}{2}}}.
\end{cases}
\end{align*}
Now, we bound I, II and III. Note that I  and III enjoy the same upper bound as
\begin{align*}
\textrm{I}(\textrm{III})\leq  2\EE\Big[\|\ttheta^{k}-\ttheta^{*} \|\cdot\|  \ttheta^k- \ttheta^{k-K_0}\|/( v^{k-1}+\delta)^{\frac{1}{2}}\Big].
\end{align*}
With Lemma \ref{legeo}, we have
$
\textrm{II}\leq \frac{2R\kappa \rho^{K_0}}{\sqrt{\delta}}.
$
Thus, we get
\begin{align}\label{yu-le15}
&\textrm{I}+\textrm{II}+\textrm{III} \nonumber\\
&\leq 4\sum_{h=1}^{K_0} \frac{\EE\Big[\|\ttheta^{k}-\ttheta^{*}\|\cdot\|\ttheta^{k+1-h}-\ttheta^{k-h}\|\Big]}{(v^{k-1} +\delta)^{\frac{1}{2}}} +\frac{2R\kappa \rho^{K_0}}{\sqrt{\delta}} \nonumber\\
& \leq 4\sum_{h=1}^{K_0} \frac{\EE\Big[\|\ttheta^{k}-\ttheta^{*}\|\cdot\|{\bf m}^{k-h}\|\Big]}{(v^{k-1} +\delta)^{\frac{1}{2}}(v^{k-h}+\delta)^{\frac{1}{2}}}+\frac{2R\kappa \rho^{K_0}}{\sqrt{\delta}}.
\end{align}
On the other hand, with the Cauchy-Schwarz inequality, we have
\begin{align}\label{yu-le2}
&\sum_{h=1}^{K_0}\EE\left(\frac{\|\ttheta^{k}-\ttheta^{*}\|\cdot\|{\bf m}^{k-h}\|}{(v^{k-1}+\delta)^{\frac{1}{2}}\cdot(v^{k-h}+\delta)^{\frac{1}{2}}}\mid\sigma^k\right)\nonumber\\
&\leq\frac{1}{\delta^{1/4}}\sum_{h=1}^{K_0}\frac{\|\ttheta^{k}-\ttheta^{*}\|}{(v^{k-1}+\delta)^{1/4}}\cdot\frac{\|{\bf m}^{k-h}\|}{( v^{k-h}+\delta)^{1/2}}\nonumber\\
&\leq  \frac{1}{2\delta^{1/4}}\sum_{h=1}^{K_0}\Big[\frac{\delta^{1/4}(1-\gamma)\omega}{4K_0}\frac{\|\ttheta^{k}-\ttheta^{*}\|^2}{(  v^{k-1}+\delta)^{1/2}}\nonumber\\
&\qquad+\frac{4K_0}{\delta^{1/4}(1-\gamma)\omega}\frac{\|{\bf m}^{k-h}\|^2}{(  v^{k-h} +\delta)}\Big]\nonumber\\
&=\frac{(1-\gamma)\omega}{8}\frac{\|\ttheta^{k}-\ttheta^{*}\|^2}{(v^{k-1}+\delta)^{1/2}}+\frac{2K_0\sum_{h=K_0}^1\frac{\|{\bf m}^{k-h}\|^2}{(  v^{k-h} +\delta)}}{ \delta^{1/2}(1-\gamma)\omega}\nonumber\\
&\leq\frac{1}{8}\frac{\langle \ttheta^*-\ttheta^k, {\bf g}(\ttheta^k)\rangle}{(v^{k-1}+\delta)^{1/2}}+\frac{2K_0\sum_{h=K_0}^1\frac{\|{\bf m}^{k-h}\|^2}{(  v^{k-h} +\delta)}}{\delta^{1/2}(1-\gamma)\omega}.
\end{align}
Combining \eqref{yu-le1}, \eqref{yu-le15} and \eqref{yu-le2}, we then get the result.

\subsection{Proof of Lemma \ref{core0}}

For simplicity, we define ${\bf g}^k:=\overline{{\bf g}}(\ttheta^k;s_{k},s_{k+1})$.
Recall ${\bf m}^k=(1-\beta)\sum_{j=1}^{k-1}\beta^{k-1-j}{\bf g}^j.$
We have
\begin{align*}
&\|{\bf m}^{k}\|^2/( v^{k}+\delta)\\
& \leq  (1-\beta)^2 \|\sum_{j=1}^{k-1}\beta^{k-1-j}{\bf g}^j /(  v^{k} +\delta)^{\frac{1}{2}}\|^2\\
&\overset{a)}{\leq}   (1-\beta)^2 (\sum_{j=1}^{k-1}\beta^{k-1-j})\cdot\sum_{j=1}^{k-1}\beta^{k-1-j}\|{\bf g}^j\|^2/( v^{k} +\delta)\\
&=(1-\beta)\cdot\sum_{j=1}^{k-1}\beta^{k-1-j}\|{\bf g}^j\|^2/( v^{k}+\delta) \\
&\overset{b)}{=}(1-\beta)\cdot\sum_{j=1}^{k-1}\beta^{k-1-j}\|{\bf g}^j\|^2/( v^{j}+\delta)
\end{align*}
where $a)$ uses the fact
$(\sum_{j=1}^{k-1}a_j b_{j})^2\leq \sum_{j=1}^{k-1}a_j^2 \sum_{j=1}^{k-1} b_{j}^2$ with $a_j=\beta^{\frac{k-1-j}{2}}$ and $b_{j}=\beta^{\frac{k-1-j}{2}}{\bf g}^j/( v^k +\delta)^{\frac{1}{2}}$, and $b)$ is due to $ v^{j}\leq  v^{k}$ when $j\leq k-1$. Then, we then get
\begin{align*}
&\sum_{k=1}^{K}\sum_{j=1}^{k-1}\beta^{k-1-j}\|{\bf g}^j\|^2/( v^{j}+\delta)\\
&=\sum_{j=1}^{K-1}\sum_{k=j}^{K-1}\beta^{k-j}\|{\bf g}^j\|^2/( v^{j}+\delta)
\leq \sum_{j=1}^{K-1}\|{\bf g}^j\|^2/(v^{j}+\delta).
\end{align*}
Combining the inequalities above, we then get the result. By applying Lemma \ref{le1}, we then get the second bound.

\subsection{Proof of Lemma \ref{core1}}
  With direct computations, we have
\begin{align*}
\Delta_k
&=\underbrace{\EE\left( \langle \ttheta^{k}-\ttheta^{*},  {\bf m}^k\rangle/( v^{k-1}+\delta)^{\frac{1}{2}}\right)}_{\textrm{I}}\\
&\quad+\underbrace{\EE\left( \frac{\langle \ttheta^{k}-\ttheta^{*},  {\bf m}^k\rangle}{(v^{k}+\delta)^{\frac{1}{2}}}-\frac{\langle \ttheta^{k}-\ttheta^{*},  {\bf m}^k\rangle}{( v^{k-1}+\delta)^{\frac{1}{2}}}\right)}_{\textrm{II}}
\end{align*}
Now, we bound I and II. The Cauchy's inequality then gives us
\begin{align*}
 \textrm{II}
 &\leq \|\ttheta^{k}-\ttheta^{*}\|\cdot\|{\bf m}^k\|\cdot(1/(v^{k-1}+\delta)^{\frac{1}{2}}-1/(v^{k}+\delta)^{\frac{1}{2}})\\
 &\leq  2RG\frac{\|{\bf g}^k\|^2}{(v^{k}+\delta)^{\frac{1}{2}}(v^{k-1}+\delta)^{\frac{1}{2}}((v^{k}+\delta)^{\frac{1}{2}}+(v^{k-1}+\delta)^{\frac{1}{2}})}\\
 &=\frac{2RG}{\delta} \cdot\frac{\|{\bf g}^k\|^2}{v^k+\delta}.
\end{align*}
 With scheme of the algorithm, we
 use a shorthand notation $\Lambda:=\EE(\langle \ttheta^{k}-\ttheta^{*}, {\bf g}^k\rangle/( v^{k-1}+\delta)^{\frac{1}{2}}\mid\sigma^k)$ and then
 get
\begin{align*}
& \textrm{I}= \EE\left(\langle \ttheta^{k}-\ttheta^{*}, \beta {\bf m}^{k-1} + (1-\beta) {\bf g}^k\rangle/( v^{k-1}+\delta)^{\frac{1}{2}}\right)\\
&=(\beta-1)\Lambda+\beta\langle \ttheta^{k}-\ttheta^{*},  {\bf m}^{k-1}\rangle/( v^{k-1}+\delta)^{\frac{1}{2}}\\
&=(\beta-1)\Lambda+\beta \Delta_k+ \beta \langle \ttheta^{k}-\ttheta^{k-1},  {\bf m}^{k-1}\rangle/( v^{k-1}+\delta)^{\frac{1}{2}}\\
&\overset{a)}{\leq}(\beta-1)\Lambda+\beta \Delta_k+\beta  \|\ttheta^{k-1}-\ttheta^{k}\|\cdot\| {\bf m}^{k-1}\|/( v^{k-1}+\delta)^{\frac{1}{2}}\\
&\overset{b)}{\leq}(\beta-1)\Lambda+\beta \Delta_k+\eta\beta  \|{\bf m}^{k-1}\|^2/(v^{k-1}+\delta),
\end{align*}

where $a)$ uses the Cauchy's inequality,  $b)$ depends on the scheme of AdaTD(0).
Combination of  the inequalities I, II and Lemma \ref{lemmayu} gives the final result.

The bound for $\Delta_k$ a direct result from the bound of ${\bf m}^k$ and ${\bf g}^k$.


\section{  Proofs for AdaTD($\lambda$)}

\subsection{Technical Lemmas}
Although $({\bf m}^k)_{k\geq 0}$ and  $(v^k)_{k\geq 0}$ are generated as the same as AdaTD(0), the difference scheme of updating $({\bf g}^k)_{k\geq 0}$ makes they different. Thus, we use different notation here.
Like previous proofs, we denote three items for AdaTD($\lambda$) as follows
\begin{align}\label{notation2}
     \begin{cases}
    \hat{\hat{{\bf m}}}_k&:=  \EE \Big[\|{\bf m}^{k}\|^2/(v^{k}+\delta)\Big],\\
    \hat{\Delta}_k&:=\EE\left(\langle \ttheta^{k}-\ttheta^{*},  {\bf m}^k\rangle/(  v^k+\delta)^{\frac{1}{2}}\right),\\
    \hat{\phi}_k&:=\EE\left(\langle \ttheta^{*}-\ttheta^{k},  {\bf g}^{\lambda}(\ttheta^k)\rangle/(  v^{k-1}+\delta)^{\frac{1}{2}}\right),\\
\hat{\Re}_k&:= \frac{8K_0(1-\beta)}{\delta^{1/2}(1-\gamma)\omega(1-\gamma\lambda)}\sum_{d=K_0}^1\hat{\hat{{\bf m}}}_{k-d}+\frac{2R(1-\beta)\kappa \rho^{K_0}}{\sqrt{\delta}(1-\gamma\lambda)}\\
&+2R\hat{G}(1-\beta)\frac{d}{\delta} \cdot\frac{\|{\bf g}^k\|^2}{v^k+\delta}
+\frac{2R(1-\beta)\hat{G} }{\sqrt{\delta}(1-\gamma\lambda)}\zeta_{k}+\eta\beta\hat{\hat{{\bf m}}}_{k}.
\end{cases}
\end{align}
Direct computing gives \begin{align}\label{summ}
\sum_{k=1}^{\infty}\zeta_k\leq \frac{\max\{\gamma\lambda,\rho\}}{(1-\max\{\gamma\lambda,\rho\})^2}.
\end{align}

\begin{lemma}\label{boundzk}
Assume $(\EE{\bf z}^k)_{k\geq 0}$ is generated by AdaTD($\lambda$). It then holds that
\begin{align}
	\|\EE {\bf z}^{k}-\EE {\bf z}^{\infty}\|\leq \frac{k\zeta^{k}}{1-\gamma\lambda}~~~{\rm and}~~~\|{\bf z}^{\infty}\|\leq \frac{1}{1-\gamma\lambda}.
\end{align}
\end{lemma}

\begin{lemma}\label{hcore0}
Assume  $({\bf m}^k)_{k\geq 0}$ and  $(v^k)_{k\geq 0}$ are given by \eqref{notation2}.
We have
\begin{align*}
\sum_{k=1}^K\hat{\hat{ {\bf m}}}_k\leq \sum_{j=1}^{K-1}\EE\|{\bf g}^j\|^2/( v ^{j}+\delta).
\end{align*}
Further, with the boundedness of $(\|{\bf g}^k\|)_{k\geq 0}$, we then get
\begin{align*}
\sum_{k=1}^K\hat{\hat{{\bf m}}}_k\leq  \ln(\frac{\delta+(K-1)\hat{G}^2}{\delta}).
\end{align*}
\end{lemma}

\begin{lemma}\label{legeo2}
Assume $(\ttheta^k)_{k\geq 0}$ is generated by AdaTD($\lambda$). Given integer $K_0\in \mathbb{Z}^+$, we then get
\begin{align*}
&\Big\|\EE\Big[ \overline{{\bf g}}^\lambda(\ttheta^{k-K_0};s_{k},s_{k+1},{\bf z}^k)\big]- {\bf g}^{\lambda}(\ttheta^{k-K_0})\Big\|\\
&\quad\leq  \frac{\kappa}{1-\gamma\lambda} \rho^{K_0}+\frac{\hat{G}}{1-\gamma\lambda} \zeta_k.
\end{align*}

\end{lemma}

\begin{lemma}\label{hlemmayu}
Assume  $(\ttheta^k)_{k\geq 0}$ is generated by AdaTD($\lambda$). Given $K_0\in \mathbb{Z}^+$, we have
\begin{align}
&\EE\Big[\frac{\langle  \ttheta^{k}-\ttheta^{*},  {\bf g}^k\rangle}{( v^{k-1}+\delta)^{\frac{1}{2}}}\Big]\leq-\frac{1}{2}\hat{\phi}_k+\frac{8K_0}{\delta^{1/2}(1-\gamma)\omega(1-\gamma\lambda)}\nonumber\\
&+\frac{2R }{\sqrt{\delta}}\Big( \frac{\kappa}{1-\gamma\lambda} \rho^{K_0}+\frac{\hat{G}}{1-\gamma\lambda} \zeta_{k}\Big).
\end{align}
\end{lemma}

\begin{lemma}\label{hcore1}
Assume $   \hat{\hat{\bf m}}_k, \hat{\Delta}_k$ are denoted by \eqref{notation2},
we then have the following result
\begin{align}
\hat{\Delta}_k+
 \frac{(1-\beta)}{2}\hat{\phi}_k\leq\beta\hat{\Delta}_{k-1}+\hat{\Re}_k.
\end{align}
On the other hand, we have another bound for $\hat{\Delta}_k$ as
$|\hat{\Delta}_k|\leq \eta\frac{2R\hat{G}}{\delta^{\frac{1}{2}}}.$
\end{lemma}

\subsection{Proof of Theorem \ref{th2}}
Similar to \eqref{lat}, we can derive
\begin{align}\label{lamdaacc}
&\sum_{k=K_0+1}^K \frac{\langle\ttheta^{*}-\ttheta^{k}, {\bf g}^{\lambda}(\ttheta^{k})\rangle}{( v^{k-1}+\delta)^{\frac{1}{2}}}\leq \frac{\EE\|\ttheta^*-\ttheta^{K_0}\|^2}{\eta}\nonumber\\
&\quad+\eta\sum_{k=K_0}^{K-1}\hat{\hat{{\bf m}}}_k+\frac{2}{1-\beta}\sum_{k=K_0+1}^K\hat{\Re}_k+\frac{4\hat{R}\hat{G}\eta}{\delta^{\frac{1}{2}}(1-\beta)}.
\end{align}
 Lemma \ref{hcore0} give us
\begin{align*}
&\eta\sum_{k=K_0}^{K-1}\hat{\hat{{\bf m}}}_k+\frac{2}{1-\beta}\sum_{k=K_0+1}^K\hat{\Re}_k\nonumber\\
&\leq\eta\sum_{k=K_0}^{K-1}\hat{\hat{{\bf m}}}_k+ \left(\frac{16K_0^2}{\delta^{1/2}(1-\gamma)\omega}+\frac{2\eta\beta}{1-\beta}\right)\sum_{k=K_0+1}^K\hat{\hat{{\bf m}}}_k\nonumber\\
&\quad+\frac{4R\hat{G}}{\delta} \sum_{k=K_0+1}^K\frac{\|{\bf g}^k\|^2}{v^k+\delta}+\frac{4R \kappa \rho^{K_0}(K-K_0)}{(1-\gamma\lambda)\sqrt{\delta}}\nonumber\\
&\quad+\sum_{k=K_0}^{K-1}\frac{2\hat{R}\hat{G} }{\sqrt{\delta}(1-\gamma\lambda)}\zeta_{k} \nonumber\\
\end{align*}
Further with Lemma \ref{le1}, the right hand is bounded by
\begin{align*}
&\Big( \frac{16K_0^2}{\delta^{1/2}(1-\gamma)\omega(1-\gamma\lambda)} +\frac{2\eta\beta}{(1-\beta)}+\eta\Big)\times\ln(\frac{\delta+(K-1)\hat{G}^2}{\delta})\nonumber\\
&+\frac{4\hat{R}\hat{G}}{\delta} \ln(\frac{\delta+K\hat{G}^2}{\delta})+\frac{4\hat{R} \kappa \rho^{K_0}(K-K_0)}{\sqrt{\delta}(1-\gamma\lambda)}+\frac{2\hat{R}\hat{G}\sum_{k=1}^{+\infty}\zeta_{k} }{\sqrt{\delta}(1-\gamma\lambda)}.
\end{align*}
The rest of the proof is almost identical to the one of Theorem \ref{th1}.
Letting $K_0=\ulcorner\ln K/\ln\frac{1}{\rho}\urcorner$ and
\begin{align*}
 &C_1^{\lambda}:=\frac{16K_0^2\hat{G}}{\delta^{1/2}(1-\gamma)(1-\alpha)\omega^2}+\frac{2\eta\beta\hat{G}}{(1-\beta)(1-\alpha)\omega}\\
 &\quad+\frac{\eta\hat{G}}{(1-\alpha)\omega}+\frac{4R\hat{G}^2}{\delta(1-\alpha)\omega}=O(\frac{(\ln K/\ln\frac{1}{\rho})^2}{\sqrt{\delta}(1-\gamma\lambda)}),\\
 &C_2^{\lambda}:=\frac{4\hat{R}^2\hat{G}}{\eta(1-\alpha)\omega}+\frac{4\hat{R}\hat{G}^2\eta}{\delta^{\frac{1}{2}}(1-\beta)(1-\alpha)\omega}\\
 &\quad+\frac{2\hat{R}\hat{G}^2\sum_{k=1}^{+\infty}\zeta_{k}}{\sqrt{\delta}(1-\gamma\lambda)(1-\alpha)\omega}+\frac{4\hat{R} \bar{\kappa}(B+\hat{R}\gamma+\hat{R})\hat{G} }{(1-\alpha)\sqrt{\delta}\omega}\\
 &\quad=O(\frac{1}{\sqrt{\delta}(1-\beta)}+\frac{\max\{\gamma\lambda,\rho\}  }{\sqrt{\delta}(1-\gamma\lambda)(1-\max\{\gamma\lambda,\rho\})^2}).
\end{align*}
we then proved the results.

\subsection{Proof of Proposition \ref{pro1}}
Based on Theorem 3, in the general case, we directly get the complexity of AdaTD($\lambda$).
Beginning from \eqref{lamdaacc}, like the proof of Theorem \ref{th1a}, the improved complexity can be obtained with condition \eqref{fastass}.

\section{Proofs of Technical Lemmas for Ada-TD($\lambda$)}
\subsection{Proof of Lemma \ref{boundzk}}
With the scheme of updating ${\bf z}^k$, it follows
$\EE {\bf z}^{k}=\sum_{t=1}^{k}(\gamma\lambda)^{k-t}\EE\phi(s_{t}).$
Assume $\pi^0$ is the initial probability  of $s_0$, then $\EE\phi(s_{t})=\Phi^{\top}\Pc^t\pi^0$ which yields
$$\EE {\bf z}^{k}=\sum_{t=1}^{k}(\gamma\lambda)^{k-t}\Phi^{\top}\Pc^t\pi^0.$$
From \cite{oldenburger1940infinite}, there exists orthogonal matrix $U$ such that
$$\Pc=U^{\top}\textrm{Diag}(1,\lambda_2,\ldots, \lambda_S)U.$$
with $\lambda_2\geq\ldots\geq\lambda_S\geq0$ and $\lambda_2\leq \rho$. Without loss of generation, we assume $\gamma\lambda>\rho$  and then derive
$$\EE {\bf z}^{k}=\Phi^{\top}U^{\top}\Gamma U\pi^0,$$ where $$\Gamma:=\textrm{Diag}(\sum_{t=1}^{k}(\gamma\lambda)^{k-t},\sum_{t=1}^{k}(\gamma\lambda)^{k-t}\lambda_2^t,\ldots,\sum_{t=1}^{k}(\gamma\lambda)^{k-t}\lambda_S^t).$$
We also see
$\EE {\bf z}^{\infty}=\Phi^{\top}U^{\top}\textrm{Diag}( \frac{1}{1-\gamma\lambda},0,\ldots, 0)U\pi^0.$
And hence, it follows
$$\EE {\bf z}^{k}-\EE {\bf z}^{\infty}=\Phi^{\top}U^{\top}
\hat{\Gamma}
U\pi^0,$$
where
$$\hat{\Gamma}:=\textrm{Diag}( \frac{-(\gamma\lambda)^k}{1-\gamma\lambda}, \sum_{t=1}^{k}(\gamma\lambda)^{k-t}\lambda_2^t,\ldots, \sum_{t=1}^{k}(\gamma\lambda)^{k-t}\lambda_S^t).$$
For $i\in\{2,3,\ldots,S\}$,
$\sum_{t=1}^{k}(\gamma\lambda)^{k-t}\lambda_S^t\leq \sum_{t=1}^{k}(\gamma\lambda)^{k-t}\rho^t.$
It is easy to see
\begin{align*}
   \zeta_k&:=\sum_{t=1}^{k}(\gamma\lambda)^{k-t}\rho^t
\leq \sum_{t=1}^k  (\max\{\gamma\lambda,\rho\})^{k-t}(\max\{\gamma\lambda,\rho\})^{t}\\
&=k(\max\{\gamma\lambda,\rho\})^{k}.
\end{align*}

\subsection{Proof of Lemma \ref{hcore0}}
This proof is identical to the one of Lemma \ref{core0} and will not be reproduced.

\subsection{Proof of Lemma \ref{legeo2}}
With the boundedness, we are led to
\begin{align*}
&\Big\|\EE \Big[ \overline{{\bf g}}^\lambda(\ttheta^{k-K_0};s_{k},s_{k+1},{\bf z}^k)\Big] \\
&-\EE  \Big[ {\bf g}^{\lambda}(\ttheta^{k-K_0};s_{k},s_{k+1},{\bf z}^{\infty})\Big]\Big\|\leq \frac{\hat{G}}{1-\gamma\lambda}k\zeta_{k}.
\end{align*}
By using a shorthand notation $$\Upsilon:=\Rc(s,s')+\gamma \phi(s')^\top \ttheta^{k-K_0}- \phi(s)^\top \ttheta^{k-K_0},$$ with direct computation, we have
\begin{align}\label{biasedexp2}
&\EE\Big[\overline{{\bf g}}^\lambda(\ttheta^{k-K_0};s_{k},s_{k+1},{\bf z}^{\infty})\mid \sigma(\sigma^{k-K_0},{\bf z}^{\infty})\Big]\nonumber\\
& =\sum_{s, s' \in \Sc} \pi(s) \Pc(s'|s) \cdot\Upsilon\cdot{\bf z}^{\infty} \nonumber\\
& +\sum_{s, s' \in \Sc} (\Pc(s_k\mid s_{k-K_0}=s)-\pi(s)) \Pc(s'|s) \cdot\Upsilon\cdot{\bf z}^{\infty} .
\end{align}
Noticing that
\begin{align*}
    &\sum_{s, s' \in \Sc} \pi(s) \Pc(s'|s) \cdot\Upsilon\cdot{\bf z}^{\infty}\\
    &=\EE\Big[{\bf g}^{\lambda}(\ttheta^{k-K_0};s_{k},s_{k+1},{\bf z}^{\infty})\mid \sigma(\sigma^{k-K_0},{\bf z}^{\infty})\Big].
\end{align*}
Direct calculation gives us
\begin{align*}
&\Big\|\sum_{s, s' \in \Sc} (\Pc(s_k\mid s_{k-K_0}=s)-\pi(s)) \Pc(s'|s)\cdot\Upsilon\cdot{\bf z}^{\infty}\Big\|\\ &\leq\frac{\kappa\rho^{K_0}}{1-\gamma\lambda}.
\end{align*}

\subsection{Proof of Lemma \ref{hlemmayu}}
Noting
$$\EE\Big[ \overline{{\bf g}}^\lambda(\ttheta;\hat{s}_{k},\hat{s}_{k+1},{\bf z}^{\infty})\Big]={\bf g}^{\lambda}(\ttheta), $$
where $(\hat{s}_{1},\hat{s}_{2},\ldots)$ is the stationary sequence. Then for any $\ttheta^1,\ttheta^2\in\RR^d$,
\begin{align*}
&\|{\bf g}^{\lambda}(\ttheta^1)-{\bf g}^{\lambda}(\ttheta^2)\|\\
&=\|\EE \Big[\overline{{\bf g}}^\lambda(\ttheta^1;\hat{s}_{k},\hat{s}_{k+1},{\bf z}^{\infty})\Big]-\EE \Big[\overline{{\bf g}}^\lambda(\ttheta^2;\hat{s}_{k},\hat{s}_{k+1},{\bf z}^{\infty})\Big]\|\\
&\leq \frac{2\|\ttheta^1-\ttheta^2\|}{1-\gamma\lambda}.
\end{align*}
With the fact $ {\bf g}^k={\bf g}^{\lambda}(\ttheta^{k})+{\bf g}^k-\overline{{\bf g}}^\lambda(\ttheta^{k-K_0};s_{k},s_{k+1},{\bf z}^k) +\overline{{\bf g}}^\lambda(\ttheta^{k-K_0};s_{k},s_{k+1},{\bf z}^k)-{\bf g}^{\lambda}(\ttheta^{k-K_0}) +{\bf g}^{\lambda}(\ttheta^{k-K_0})-{\bf g}^{\lambda}(\ttheta^{k})$, we can have
\begin{align}\label{hyu-le1}
&\EE\Big[\langle \ttheta^{k}-\ttheta^{*},  {\bf g}^k\rangle/( v^{k-1}+\delta)^{\frac{1}{2}}\Big]\nonumber\\
&=\EE\Big[\langle\ttheta^{k}-\ttheta^{*}, {\bf g}^{\lambda}(\ttheta^{k})\rangle/(v^{k-1}+\delta)^{\frac{1}{2}}\Big]\nonumber\\
&+\underbrace{\frac{\EE\mid\langle \ttheta^{k}-\ttheta^{*},  \left[ {\bf g}^k-\overline{{\bf g}}^\lambda(\ttheta^{k-K_0};s_{k},s_{k+1},{\bf z}^k)\right]\rangle\mid}{(v^{k-1}+\delta)^{\frac{1}{2}}}}_{\textrm{I}}\nonumber\\
&+\underbrace{\frac{\EE\mid\langle \ttheta^{k}-\ttheta^{*},  \left[ \overline{{\bf g}}^\lambda(\ttheta^{k-K_0};s_{k},s_{k+1},{\bf z}^k)-{\bf g}^{\lambda}(\ttheta^{k-K_0})\right]\rangle\mid}{(v^{k-1}+\delta)^{\frac{1}{2}}}}_{\textrm{II}}\nonumber\\
&+\underbrace{\frac{\EE\mid\langle \ttheta^{k}-\ttheta^{*},  \left[ {\bf g}^{\lambda}(\ttheta^{k-K_0})-{\bf g}^{\lambda}(\ttheta^{k})\right]\rangle\mid}{(v^{k-1}+\delta)^{\frac{1}{2}}}}_{\textrm{III}}.
\end{align}
We  bound I, II and III in the following. We can see I and III have the same bound
\begin{align*}
&\textrm{I}(\textrm{III})\leq  \frac{2}{1-\gamma\lambda}\EE\Big[\frac{\|\ttheta^{k}-\ttheta^{*} \|\cdot\|  \ttheta^k- \ttheta^{k-K_0}\|}{( v^{k-1}+\delta)^{\frac{1}{2}}}\Big],
\end{align*}
and with Lemma \ref{legeo2}
\begin{align*}
\textrm{II}\leq \frac{2\hat{R} }{\sqrt{\delta}}( \frac{\kappa}{1-\gamma\lambda} \rho^{K_0}+\frac{\hat{G}}{1-\gamma\lambda} \zeta_{k}).
\end{align*}
Hence, we have
\begin{align}
&\textrm{I}+\textrm{II}+\textrm{III} \nonumber\\
&\leq \frac{4}{1-\gamma\lambda}\sum_{h=K_0}^1 \frac{\EE\Big[\|\ttheta^{k}-\ttheta^{*}\|\cdot\|\ttheta^{k+1-h}-\ttheta^{k-h}\|\Big]}{(v^{k-1} +\delta)^{\frac{1}{2}}} \nonumber\\
&\quad+\frac{2\hat{R} }{\sqrt{\delta}}\Big( \frac{\kappa}{1-\gamma\lambda} \rho^{K_0}+\frac{\hat{G}}{1-\gamma\lambda} \zeta_{k}\Big) \nonumber\\
&\leq \frac{4}{1-\gamma\lambda}\sum_{h=K_0}^1 \frac{\EE\Big[\|\ttheta^{k}-\ttheta^{*}\|\cdot\|{\bf m}^{k-h}\|\Big]}{(v^{k-1} +\delta)^{\frac{1}{2}}\cdot(v^{k-h}+\delta)^{\frac{1}{2}}}\nonumber\\
&\quad+\frac{2\hat{R} }{\sqrt{\delta}}\Big( \frac{\kappa}{1-\gamma\lambda} \rho^{K_0}+\frac{\hat{G}}{1-\gamma\lambda} \zeta_{k}\Big).
\end{align}
On the other hand, with the Cauchy-Schwarz inequality, we derive
\begin{align}\label{hyu-le2pre}
&\sum_{h=K_0}^1\EE\left(\frac{\|\ttheta^{k}-\ttheta^{*}\|\cdot\|{\bf m}^{k-h}\|}{(v^{k-1}+\delta)^{\frac{1}{2}}\cdot(v^{k-h}+\delta)^{\frac{1}{2}}}\mid\sigma^k\right)\nonumber\\
&\leq\frac{1}{\delta^{1/4}}\sum_{h=K_0}^1\frac{\|\ttheta^{k}-\ttheta^{*}\|}{(v^{k-1}+\delta)^{1/4}}\cdot\frac{\|{\bf m}^{k-h}\|}{( v^{k-h}+\delta)^{1/2}}\nonumber\\
&\leq  \frac{1}{2\delta^{1/4}}\sum_{h=K_0}^1\Big(\frac{\delta^{1/4}(1-\alpha)\omega(1-\gamma\lambda)}{4K_0}\frac{\|\ttheta^{k}-\ttheta^{*}\|^2}{(  v^{k-1}+\delta)^{1/2}}\nonumber\\
&\quad+\frac{4K_0}{\delta^{1/4}(1-\alpha)\omega(1-\gamma\lambda)}\frac{\|{\bf m}^{k-h}\|^2}{(  v^{k-h} +\delta)}\Big).
\end{align}

The right hand of \eqref{hyu-le2pre} is further bounded by
\begin{align}\label{hyu-le2}
&\frac{(1-\alpha)\omega(1-\gamma\lambda)}{8}\|\ttheta^{k}-\ttheta^{*}\|^2/(v^{k-1}+\delta)^{1/2}\nonumber\\
&\quad+\frac{2K_0}{\delta^{1/2}(1-\alpha)\omega(1-\gamma\lambda)}\sum_{h=1}^{K_0} \hat{\hat{{\bf m}}}_{k-h}\nonumber\\
&\leq\frac{(1-\gamma\lambda)}{8}\hat{\phi}_k +\frac{2K_0}{\delta^{1/2}(1-\alpha)\omega(1-\gamma\lambda)}\sum_{h=1}^{K_0} \hat{\hat{{\bf m}}}_{k-h}.
\end{align}

Combining \eqref{hyu-le1} and \eqref{hyu-le2}, we then get the result.

\subsection{Proof of Lemma \ref{hcore1}}
The proof is identical to Lemma \ref{core1} and will not be repeated.
\end{document}